# Searching for Line Transitive, Point Imprimitive, Linear Spaces [1]

Gregory Cresp

October 2001

[1] Submitted in partial fulfilment of the requirements for the Bachelor of Science degree with Honours at The University of Western Australia

## Abstract

A finite linear space is a finite set of points and lines, where any two points lie on a unique line. Well known examples include projective planes. This project focuses on linear spaces which admit certain types of symmetries. Symmetries of the space which preserve the line structure are called *automorphisms*. A group of these is called an *automorphism group* of the linear space. Two interesting properties of linear spaces are point imprimitivity and line transitivity. Point imprimitive spaces admit a second structure on the points aside from the lines, which is also preserved by an automorphism group. In line transitive spaces, given any two lines, an automorphism can be found that maps one line to the other.

Very few point imprimitive, line transitive linear spaces, apart from projective planes, are known. Such spaces that have been found have been surprising. One point of interest is whether such spaces are rare and the known ones are in some sense exceptional, or if there are many such spaces, but mathematicians have been looking in the wrong places.

Here we investigate methods to construct a line transitive, point imprimitive linear space over a given point set and automorphism group. We employ these methods on two given automorphism groups, both on a set of 451 points. This was an exceptional situation identified in theoretical work of Praeger and Tuan. Included in this is the development of an algorithm, written in GAP, an algebraic programming system, and C, to perform these constructions. This algorithm is extendible to a wider class of groups.

**Acknowledgements**

There are many people who contributed to the production of this dissertation.

My thanks go to my supervisors Prof. Cheryl E. Praeger and Dr. Alice Niemeyer. Without their time and effort this work would never have even begun.

I would like to also acknowledge the assistance of Dr. Anton Betten. Much of the C code used in this work is based on code written by him, and his advice on design theory was invaluable. My thanks go to Mr. Maska Law, who assisted with background on projective planes of order 7, and Dr. Rick McFeat, who directed me on Galios theory.

For their efforts in proof reading and their continued support throughout the year, my sincerest thanks for Dr. John Cresp, Mrs. Judy Cresp and Miss Kate Gilders.

Finally I must acknowledge the support of my fellow honours students for their advice on latex, GAP and not going insane during the year.

# Contents

# Chapter 1

# Introduction

## 1.1 Introduction

A finite linear space is a finite set of points and lines such that any two points lie on a unique line. There are two properties of particular interest of linear spaces, line transitivity and point imprimitivity, which we will introduce in more detail in Chapter 3.

We refer to maps on the space which preserve the line structure as automorphisms. A linear space is line transitive if it is so symmetric that, given any two lines, an automorphism exists which maps one of these lines to the other. If there is a non-trivial partition on points that is preserved by all automorphisms, then the space is called point imprimitive. Transitive and imprimitive are defined more rigorously in Chapter 2.

There is much interest in finding particular examples of point imprimitive, line transitive finite linear spaces. There are many projective planes (a particular class of linear spaces where any two lines meet at a point) known to have these properties. We will consider an infinite family of such planes in Section 3.5. However, there are very few known examples of spaces having these properties that are not projective planes. One point of interest here is whether those examples that are known are in some way exceptional, or whether there are many examples, but they have yet to be identified. In order to answer this question, more information and more examples need to be found.

Designs are a class of algebraic objects satisfying very particular requirements. Some designs are also linear spaces and these designs form an easily identifiable subclass of all designs, which we will see in Chapter 4. One way of finding linear spaces of the type required is to construct designs which are linear spaces and satisfy the required properties.

One way to construct designs uses a result known as the Orbit Theorem (Theorem 4.3.1). This theorem shows that a design can be constructed from a group action on a set if a subset of a particular size can be found with particular properties. One approach to this is to write a program to search for such a subset. The group, set and some design parameters are the initial conditions of the search. Usually, there are so many subsets of the desired size that to check all of them would be impossible. Group theory, design theory and combinatorics are used before and during the search to narrow down the search space. Of particular interest is the Delandtsheer-Doyen Theorem (Theorem 4.4.2), which gives further conditions on the subsets, relating to any invariant partition on the point set. We will introduce how these results may be used in Chapter 4.

These approaches will be put into practise in Chapter 5. We will introduce two type of group actions on 451 points that were used in searches for designs. These searches were completed, but unfortunately found no examples. However, they serve as good examples of how the theory may be applied, and eliminate these actions from further consideration. The algorithm used to perform the second search is discussed in Chapter 6, where we see that the same algorithm may be used for a wider class of searches.

# Chapter 2

# Preliminaries

## Introduction

Our work will require a basic background in the theory of group actions and permutations. We assume the reader already has a background in group theory. Here we will introduce the idea of a group acting on a set and how this action can alternatively be viewed as a group of permutations on the set. This material will form the basis for work in later chapters. Throughout this chapter we will let $\Omega$ denote a finite set and $G$ a group. Also, $1_G$ will be used to denote the identity of $G$, or 1 if the group is clear from the context. Although most of the theory given here also applies to infinite sets and groups, we will restrict ourselves to finite sets and groups.

## 2.1 Permutation Groups

The study of permutation groups has been an important part of group theory. In fact, it can easily be shown that any group is isomorphic to some permutation group. We are interested in permutation groups on points of our linear spaces. A permutation group is a group of permutations on a given set, $\Omega$, using functional composition.

**Definition 2.1.1.** Let $\Omega$ be a set. A *permutation* of $\Omega$ is a 1-1 and onto map $\Omega \to \Omega$.

We call the group of all permutations on $\Omega$ the *symmetric group* of $\Omega$ and denote it $\mathrm{Sym}(\Omega)$. Any subgroup of $\mathrm{Sym}(\Omega)$ is called a *permutation group on* $\Omega$.

## 2.2 Group Actions

We will often have the elements of $G$ moving the elements of $\Omega$ around in some way. An example of this is how the group of rotations of the square affects the set of the vertices of the square. This movement is called a group action.

**Definition 2.2.1.** Let $G$ be a group and $\Omega$ a set. A *group action* of $G$ on $\Omega$ is an assignment $\Omega \times G \to \Omega$, satisfying:

- $(\omega, 1) \mapsto \omega$ for all $\omega \in \Omega$.
- $((\omega, g), h) \mapsto (\omega, gh)$ for all $\omega \in \Omega,\ g, h \in G$.

The image of $(\omega, g)$ is denoted $\omega^g$.

It is easy to see that any element of $G$ gives rise to a permutation of $\Omega$. Consider the map $\varphi_g : \Omega \to \Omega$ given by $\varphi_g(\omega) = \omega^g$ for all $\omega \in \Omega$. Clearly this is a map from $\Omega$ into itself. Since $\Omega$ is a finite set, $\varphi_g$ is 1-1 if and only if it is onto. Suppose $\varphi_g(\alpha) = \varphi_g(\beta)$. Then $\alpha^g = \beta^g$, so $(\alpha^g)^{g^{-1}} = (\beta^g)^{g^{-1}}$ and thus $\alpha = \beta$. We can extend this idea to a function $\varphi : G \to \mathrm{Sym}(\Omega)$ by $\varphi(g) = \varphi_g$. This is well defined, since each $\varphi_g$ is. This process gives us a permutation corresponding to each element of $G$. This is referred to as the *permutation representation* of the action of $G$. We can, and will, use the group action and its permutation representation interchangeably.

This permutation representation map $\varphi$ may not be one to one. It may be that two different elements of $G$ act in the same way. Since $\varphi$ is a homomorphism, this is equivalent to the kernel of the map being non-trivial. By this we mean that there is some non-identity $g \in G$ such that $\varphi(g) = 1_{Sym(\Omega)}$. If $\varphi$ is 1-1, the action is called faithful. If an action is faithful then we can identify $G$ with $\varphi(G) \leq \mathrm{Sym}(\Omega)$.

**Definition 2.2.2.** Let $G$ act on $\Omega$. The action is said to be *faithful* if the only element of $G$ to fix all elements of $\Omega$ is 1.

Given $G$ acting on $\Omega$, we can also consider $G$ acting on the subsets of $\Omega$. We define the *set-wise image* of a subset $S \subset \Omega$ under $g \in G$ to be $S^g = \{\omega^g | \omega \in S\}$.

We are often interested in where in the set a given element can be mapped to under the action. We call the set of all possible images of a point the orbit of that point.

**Definition 2.2.3.** Given $G$ acting on $\Omega$ and some $\omega \in \Omega$, the *orbit* of $\omega$ under $G$, denoted $\omega^G$ is the set $\{\omega^g | g \in G\}$.

It is also often interesting to consider, for a given $\omega \in \Omega$, which elements of $G$ fix $\omega$. The set of all such elements is called the *stabiliser of* $\omega$ and is denoted $G_\omega$.

## 2.3 Properties of Actions

The first property of group actions we consider is transitivity. A transitive action is one where it is possible to move from any element in the set to any other element of the set. Another way of looking at this is that $\Omega$ is one $G$-orbit.

**Definition 2.3.1.** Let $G$ be a group acting on $\Omega$. The action is *transitive* if, for any $\alpha, \beta \in \Omega$ there is some $g \in G$ such that $\alpha^g = \beta$.

The second property is primitivity. Before we can describe a primitive action, we need to consider partitions and invariant partitions.

**Definition 2.3.2.** A *partition* of a set $\Omega$ is a set of pairwise disjoint subsets of $\Omega$, say $\mathcal{C} = \{C_1, C_2, \cdots, C_n\}$ such that $\cup_{i=1}^{n} C_i = \Omega$.

We call the subsets in a partition classes. An alternative name is blocks, but we will see in Chapter 4 that this would cause confusion with lines. A partition is said to be $G$-invariant invariant partitions. A transitive action that is not primitive is called *imprimitive*. We will see later that we are interested in imprimitive actions.

In many applications, we wish to find the size of an orbit without having to explicitly calculate that orbit. One way of doing this is to find the stabiliser of an element in the orbit.

**Theorem 2.3.3 (Orbit Stabiliser Theorem).** *Let $G$ act faithfully on $\Omega$, a finite set and $\omega \in \Omega$. Then $|\omega^G| = \frac{|G|}{|G_\omega|}$.*

For a proof of this see [4, Theorem 1.4A].

A property of interest, but less importance to us, is regularity. We say an action is *semi-regular* if the point stabiliser of any point fixes all points. A *regular* action is one that is semi-regular and transitive.

**Lemma 2.3.4.** *Let $G$ act faitfully and regularly on $\Omega$. Then we can identify $\Omega$ with $G$ in such a way that $G$ acts by its group operation.*

**Proof.** Since $G$ is faithful and regular, the only element of $G$ that can fix any point is 1. Taking any $\omega \in \Omega$ and define the map $\phi : \Omega \to G$ by $\omega^g \mapsto g$. Since the action is transitive, any $\alpha \in \Omega$ is $\omega^g$ for some $g \in G$. Furthermore, if $\omega^g = \omega^h$ then $\omega^{gh^{-1}} = \omega$ and hence, since the point stabiliser of any point is $\{1\}$, $g = h$. Hence this map is well defined. Clearly if $\phi(\alpha) = \phi(\beta)$ then $\alpha = \omega^{\phi(\alpha)} = \omega^{\phi(\beta)} = \beta$ so the map is 1-1. Also, by the Theorem 2.3.3, $|G| = |\omega|$ and both are finite, so the map is also onto. Hence $\phi$ is an identification of $\Omega$ with $G$. Finally, notice $\phi(\alpha^g) = \phi(\alpha)g$ for all $\alpha \in \Omega, g \in G$ by the definition of $\phi$. $\square$

## Conclusion

We have introduced here the idea of a group action. There are several properties of group actions that will be important to us in later chapters. First, we will only be considering faithful actions, where the only element of the group to fix all elements of the set is 1. We have also introduced imprimitive actions and transitive actions, which will be used in defining properties of linear spaces in Chapter 3.

# Chapter 3

# Linear Spaces

## Introduction

Before considering the main topic of this thesis - searching for linear spaces with particular properties, we must first consider what constitutes a finite linear space, and how we characterise these properties. Here we introduce a linear space and define the two properties of interest: line transitivity and point imprimitivity. We will also consider the known examples of spaces with these properties. We will briefly consider projective planes, but only to explain why we decided not to search for them.

## 3.1 Finite Linear Spaces

A finite linear space is an ordered pair of sets, a set of points and a set of lines.

**Definition 3.1.1.** Let $\mathcal{P}$ be a finite set and $\mathcal{L}$ a set of subsets of $\mathcal{P}$. Elements of $\mathcal{P}$ are called *points*, elements of $\mathcal{L}$ *lines*. Then $\mathcal{S} = (\mathcal{P}, \mathcal{L})$ is a *finite linear space* if each line contains at least 2 points and, given any two points, there is a unique line containing them both.

In this dissertation we are only interested in finite linear spaces. We will always take linear space to mean finite linear space. We are mostly interested in looking at maps on the points which preserve the line structure. We call these automorphisms.

**Definition 3.1.2.** Let $\mathcal{S} = (\mathcal{P}, \mathcal{L})$ be a linear space and $g \in \mathrm{Sym}(\mathcal{P})$. Then $g$ is called an *automorphism* of $\mathcal{S}$ if, given any line $L \in \mathcal{L}$, the set-wise image, $L^g = \{p^g | p \in L\}$ is also a line.

It is clear that the composition of two automorphisms will be an automorphism, since lines will be preserved throughout. It thus makes sense to look at groups of automorphisms, using functional composition. These will be subgroups of $\mathrm{Sym}(\mathcal{P})$.

**Definition 3.1.3.** Let $\mathcal{S} = (\mathcal{P}, \mathcal{L})$ be a linear space and $G \leq \mathrm{Sym}(\mathcal{P})$. If all elements of $G$ are automorphisms of $\mathcal{S}$ then $G$ is called an *automorphism group* of $\mathcal{S}$. The group consisting of all automorphisms of $\mathcal{S}$ is called the *full automorphism group* of $\mathcal{S}$ and is denoted $\mathrm{Aut}(\mathcal{S})$.

As we will see in Chapter 4, we are often interested in linear spaces with a constant line size. When this is the case, we can construct some relationships between the number of lines and points, the number of lines per point and the number of points per line.

**Lemma 3.1.4.** *Let $\mathcal{S} = (\mathcal{P}, \mathcal{L})$ be a linear space with $|\mathcal{P}| = v$ and a constant line size $k$. Let $b = |\mathcal{L}|$ and $r$ be the number of lines through any point. Then $b = \frac{v(v-1)}{k(k-1)}$, $bk = vr$ and $v - 1 = r(k - 1)$.*

**Proof.** Consider all pairs $(\{p_1, p_2\}, L)$ where $L \in \mathcal{L}$ and $p_1, p_2 \in L$. Suppose there are $n$ of these. If we first select a pair of points, then a line containing those points we see $n = \frac{v(v-1)}{2} \cdot 1$. If we select the line first, then a pair of points on that line, $n = b\frac{k(k-1)}{2}$ hence $b = \frac{v(v-1)}{k(k-1)}$.

Now consider the pairs $(p, L)$ with $p \in \mathcal{P}, L \in \mathcal{L}, p \in L$. By first selecting $p$, then $L$ with $p \in L$ we see there are $vr$ of these. Now by selecting $L$, then $p \in L$, there are $bk$ such pairs. Hence $vr = bk$.

Combining these two equations, we see that $vr = bk = \frac{v(v-1)}{k(k-1)}k = v\frac{v-1}{k-1}$. The result follows. □

## 3.2 Projective Planes

Projective planes are a sub-class of linear spaces. In essence, a projective plane is a linear space where any two lines have a unique point in common.

**Definition 3.2.1.** Let $\mathcal{S} = (\mathcal{P}, \mathcal{L})$ be a linear space. We call $\mathcal{S}$ a *projective plane* if $|\mathcal{P}| \geq 3$, every line contains at least 3 points and any two lines intersect in a unique point.

In projective planes there is a duality between points and lines. This means that any correct statement involving points and lines is also correct if the words 'points' and 'lines' are exchanged. For example, if 'every line contains exactly $k$ points' it is also true that 'every point lies on exactly $k$ lines'.

We can see then that projective planes are very specific objects. As we will see in Section 3.5, there is an infinite family of projective planes which are linear spaces of the type we are interested in. This relates to the exceptional nature of projective planes, so in this dissertation projective planes will not be considered.

When a linear space has a constant line size, it is easy to decide whether the space is a projective plane or not. This characterisation of projective planes arises from the duality of points and lines.

**Lemma 3.2.2.** *Let $\mathcal{S} = (\mathcal{P}, \mathcal{L})$ be a linear space with constant line size. Then $\mathcal{S}$ is a projective plane if and only if $|\mathcal{P}| = |\mathcal{L}|$.*

**Proof.** Suppose $\mathcal{S}$ is a linear space with constant line size. Let $v, b, k, r$ denote the number of points, lines, points per line and lines per point respectively. Consider the tuples $(L_1, L_2, p)$ where $L_1, L_2 \in \mathcal{L}$ and $p \in L_1 \cap L_2$. Suppose there are $n$ of these. Since two lines intersect in at most one point (otherwise the linear space condition is violated), we have $n \leq b(b-1)$ and equality occurs only for projective planes. Choose $L_1$, then $p \in L_1$ then $L_2 \neq L_1$ with $p \in L_2$, so we see $n = bk(r-1)$.

First suppose $\mathcal{S}$ is a projective plane. Then we can choose $L_1$ and $L_2$ arbitrarily, they define a unique $p$. Hence $n = b(b-1) = bk(r-1)$. It follows that $b - 1 = k(r-1)$. However, $v - 1 = k(r-1)$ from Lemma 3.1.4 and hence $v = b$.

Suppose there are equal numbers of points and lines, so $v = b$. In Lemma 3.1.4 we saw $bk = vr$ and $\frac{v(v-1)}{2} = b\frac{k(k-1)}{2}$. It follows, since $v = b$, that $r = k$ and $(v-1) = k(k-1)$. Thus $n = bk(r-1) = bk(k-1) = b(v-1) = b(b-1)$, and hence $\mathcal{S}$ is a projective plane. □

## 3.3 Properties of Linear Spaces

An automorphism group of a linear space can be thought of as a group action on the points in the obvious way. We can also consider the induced action of this group on the lines, considering lines as subsets of points. Two properties often of interest in group actions are transitivity and imprimitivity. These were both discussed in Chapter 2. We will now see how these ideas lead to two important properties a linear space may have with respect to an automorphism group.

The first property that we consider is line transitivity.

**Definition 3.3.1.** Let $\mathcal{S} = (\mathcal{P}, \mathcal{L})$ be a linear space and $G \leq \mathrm{Aut}(\mathcal{S})$. Then $G$ is said to be *line transitive* on $\mathcal{S}$ if the point-wise action of $G$ on $\mathcal{L}$ is transitive. By this we mean that, given any two lines $L_1, L_2 \in \mathcal{L}$, there is some map $g \in G$ such that $L_1^g = L_2$.

We can see that if $G$ is line transitive and $G \leq H \leq \mathrm{Aut}(\mathcal{S})$ then $H$ will also be line transitive, since the required map between two lines could come from $G$. Thus line transitivity is preserved 'upwards'.

The second property of interest is point imprimitivity. We saw in Section 2.3 that a group action is said to be imprimitive if it preserves some non-trivial partition of the set.

**Definition 3.3.2.** Let $\mathcal{S} = (\mathcal{P}, \mathcal{L})$ be a linear space and $G \leq \mathrm{Aut}(\mathcal{S})$. Then $\mathcal{S}$ is said to be *$G$-point imprimitive* if there is some non-trivial $G$-invariant partition of $\mathcal{P}$.

Notice if a group leaves a partition invariant, then so will any subgroup of it. Hence if $G$ is point imprimitive on $\mathcal{S}$, then any subgroup of $G$ would also be. Thus point imprimitivity is preserved 'downwards'.

These two properties seem to be working in different directions. Line transitivity needs lines, and hence points, to move around a lot, and is preserved 'upwards'. Point imprimitivity needs points to move little under $G$ and is preserved 'downwards'. It may therefore be difficult to find some middle ground between these, where a space was be both line transitive and point imprimitive. These two properties can hold at the same time, see Section 3.5.

## 3.4 Isomorphic Linear Spaces

As with many other algebraic structures, we have a concept of isomorphism of linear spaces. Two linear spaces are isomorphic if they have the same number of points and the line structure is the same, up to some re-naming of the points.

**Definition 3.4.1.** Suppose $\mathcal{S}_1 = (\mathcal{P}, \mathcal{L})$ and $\mathcal{S}_2 = (\mathcal{Q}, \mathcal{M})$ are both linear spaces. Then $\mathcal{S}_1$ and $\mathcal{S}_2$ are said to be *isomorphic* if there is some 1-1 and onto map $\varphi : \mathcal{P} \to \mathcal{Q}$ such that $\varphi(\mathcal{L}) = \mathcal{M}$. This is denoted $\mathcal{S}_1 \cong \mathcal{S}_2$ and $\varphi$ is called an *isomorphism.*

By $\varphi(\mathcal{L}) = \mathcal{M}$ we mean that for any line $L \in \mathcal{L}$, the point-wise image of $L$ under $\varphi$ is in $\mathcal{M}$ and conversely, $\varphi^{-1}(M) \in \mathcal{L}$ for any $M \in \mathcal{M}$.

If the two linear spaces have the same point set, say $\mathcal{P}$, then the isomorphism will belong to $\mathrm{Sym}(\mathcal{P})$.

Isomorphism is an important concept, since when working with linear spaces there is no real difference between working with a space and working with another isomorphic space.

## 3.5 Examples of Point Imprimitive, Line Transitive, Linear Spaces

The first line transitive, point imprimitive, linear spaces to be identified were projective planes. In fact, an infinite family of these has been identified. This family consists of the Desarguesian projective planes under the action of Singer cycles. We will outline how such spaces and the action are constructed.

### 3.5.1 Desarguesian Planes

The non-zero elements of a field form a cyclic multiplicative group. A generator of this group is called a *primitive element* of the field.

Take some prime power $q \in \mathbb{Z}$, so $q = p^n$ for some prime $p$ and some integer $n$. Consider the field with $q$ elements, $\mathbb{F}_q$. We will use $F_q$ to construct $F_{q^3}$. Take some irreducible polynomial $f(X)$ of degree 3 in $\mathbb{F}_q[X]$. By irreducible, we mean all divisors of $f(X)$ have degree either 0 (the constant polynomials) or degree $\deg(f(X))$. Let $\alpha$ be some zero of this ploynomal. Notice $\alpha \notin \mathbb{F}_q$, since if it was, $X - \alpha$ would divide $f(X)$, contradicting irreducability. Construct the field $\mathbb{F}_{q^3}$ by $\mathbb{F}_{q^3} = \mathbb{F}_q(\alpha) = \{a + b\alpha + c\alpha | a, c, b \in \mathbb{F}_q\}$.

We say $f(X)$ is a *primitive polynomial* if $\alpha$ is a primitive element of $\mathbb{F}_{q^3}$ constructed in this way. We wish to construct $\mathbb{F}_{q^3}$ using a primitive polynomial. Find some primitive element, $\beta$, of $F_{q^3}$. Take its minimal polynomial, $f(X)$, over $F_q$ (the monic polynomial of least degree such that $f(\beta) = 0$). Identify $\beta = X$. Now construct $\mathbb{F}_{q^3}$ by taking all polynomials of degree less than 3 from $\mathbb{F}_q[X]$, with multiplication and addition modulo $f(X)$. We can think of $\mathbb{F}_{q^3}$ as $V$, a 3 dimensional vector space over $q$, with basis $\{1, X, X^2\}$. We see that this is indeed a basis, for if there were $\alpha, \beta, \gamma \in \mathbb{F}_q$ not all 0 such that $\alpha + \beta X + \gamma X^2 = 0$ (in $\mathbb{F}_{q^3}$) this would violate $f(X)$ being the minimal polynomial of $\alpha$, since $\deg(f(X))$=3. Notice from this construction, $X$ is a primivie element of $\mathbb{F}_{q^3}$.

Now define $\mathcal{S} = (\mathcal{P}, \mathcal{L})$ by $\mathcal{P}$ being the 1-dimensional subspaces of $V$ and $\mathcal{L}$ the 2-dimensional subspaces of $V$. Let a point lie on a line if the corresponding 1-D subspace is a subset of the corresponding 2-D subspace. Notice that any two distinct 1-D subspaces lie on a unique 2-D subspace and any two 2-D subspaces intersect at a unique 1-D subspace. As a result, $\mathcal{S}$ is a projective plane.

Note that any 2-D subspace contains exactly $q+1$ distinct 1-D subspaces, so the space has a constant line size. Each non-zero vector $x \in V$ lies on the unique 1-D subspace $\langle x \rangle = \{\lambda x | \lambda \in \mathbb{F}_q\}$, and each 1-D subspace contains $q - 1$ non-zero vectors. There are $q^3 - 1$ non-zero vectors total, so $|\mathcal{P}| = \frac{q^3-1}{q-1} = q^2 + q + 1$. Since $\mathcal{S}$ is a projective plane, there are also $q^2 + q + 1$ lines.

We now construct the Singer cycle, acting on $\mathcal{S}$. Take the degree 1 polynomial $X \in \mathbb{F}_{q^3}$ and consider $X$ acting on $V$ by multiplication, so $g(X) \mapsto Xg(X)$. This is an invertible linear transformation of $V$, since $X(\alpha g(X) + \beta h(X)) = \alpha X g(X) + \beta X h(X)$ by the properties of multiplication and addition on fields, and $\mathbb{F}_{q^3}$ is a field, so $(X)^{-1}$ exists. This defines an automorphism of $\mathcal{S}$, since invertible linear transformations preserve subspaces, so 1-D subspaces (points) are mapped to 1-D subspaces and 2-D subspaces (lines) to 2-D subspaces, preserving incidence.

Since $X$ is a primitive element $\mathbb{F}_{q^3}$, it is a generator of the multiplicative group $\mathbb{F}_{q^3} \setminus \{0\}$. It follows that $\langle X \rangle \cong \mathbb{Z}_{q^3-1}$ is an automorphism group of $\mathcal{S}$. Note this action on $\mathcal{P}$ is not regular, since $(q^3 - 1) > |\mathcal{P}|$.

We claim that $K$, the kernel of the action of $\langle X \rangle$ on $\mathcal{P}$ is the group of non-zero scalars of $\mathbb{F}_{q^3}$. Notice the group of non-zero scalars of $\mathbb{F}_{q^3}$ is isomorphic to $\mathbb{F}_q \setminus \{0\}$. Any element of $\mathbb{F}_q \setminus \{0\}$ fixes any point in $\mathcal{P}$. Take $p \in \mathcal{P}$, say $p = \langle x \rangle$ for some $x \in V$. Then take $\alpha \in \mathbb{F}_q \setminus \{0\}$. We see $p^\alpha := \alpha\langle x \rangle = \langle \alpha x \rangle = \langle x \rangle = p$ since $\mathbb{F}_q$ is a field. Hence $\mathbb{F}_q \setminus \{0\}$ is in the kernel of the action. We now show that any element of the kernel of $\langle X \rangle$ on $\mathcal{P}$ is also a non-zero scalar of $\mathbb{F}_{q^3}$.

Suppose we have some element of $\langle X \rangle$ that fixes all points in $\mathcal{P}$. Let this element be $X^i$ for some $i = 1, 2, \cdots, q^3 - 1$. Then, since $X^i$ fixes all points, in particular, $X^i$ must fix the subspace generated by 1, which is $\mathbb{F}_q \setminus \{0\}$. Hence $1^{X^i} \in \mathbb{F}_q \setminus \{0\}$. Notice that $X^i : 1 \mapsto X^i$ and hence $X^i \in \mathbb{F}_q \setminus \{0\}$.

As a multiplicative group, $\mathbb{F}_q \setminus \{0\} \cong \mathbb{Z}_{q-1}$. It follows that the permutation group on $\mathcal{P}$ induced by $\langle X \rangle$ is isomorphic to $\frac{\mathbb{Z}_{q^3-1}}{\mathbb{Z}_{q-1}} = \mathbb{Z}_{q^2+q+1}$. Now $|\mathcal{P}| = q^2 + q + 1$, so this action is clearly regular. Hence we can identify the points with $\mathbb{Z}_{q^2+q+1}$ and consider the action to be addition modulo $q^2 + q + 1$. The action on the lines is isomorphic to the same group and there are $q^2 + q + 1$ lines, so this action is transitive on the lines.

Suppose that $q^2 + q + 1$ is not prime, so $q^2 + q + 1 = ab$ for $a, b \in \mathbb{Z}^+$, neither being 1. We can construct a partition with $a$ classes of size $b$ that is preserved by the Singer cycle. Consider the partition $\mathcal{C} = \{C_1, C_2, \cdots, C_a\}$ where $C_i = \{bj + i | j \in \{0, 1, \cdots, a-1\}$. This partition is clearly preserved by the action since $C_i^c = C_{i+c(\text{mod } a)}$ for all $c \in \mathbb{Z}_{q^2+q+1}$. We thus have a point imprimitive, line transitive action on a linear space. There are infinitely many values of $q$ such that $q^2 + q + 1$ is not prime. Take $p \neq 3$ to be prime. There are infinitely many primes. Consider $q = p^2$. Then $p = 1$ or $2 \pmod 3$. Notice $1^2 = 2^2 = 1 \pmod 3$, so $q = 1 \pmod 3$. Hence $q^2 + q + 1 = 1 + 1 + 1 = 0 \pmod 3$, so 3 divides $q^2 + q + 1$. Clearly then $q^2 + q + 1$ is prime. Hence there are infinitely many such projective planes. However, since projective planes are already exceptional in some ways as linear spaces, the line transitive, point imprimitive spaces that are not projective planes are of more interest.

### 3.5.2 Other Linear Spaces

Of the non projective planes, only a very small number have so far been identified. The first was constructed by Mills in [6]. This was a space on 91 points, where the size of each line was 6. This space was constructed using a different method by Colbourn and Colbourn in [2] and another, non-isomorphic space on 91 points with the same line size was identified at the same time. In both cases, the construction of line transitive, point imprimitive linear spaces was not the aim of the research, it was instead a by-product of other work.

The most recently identified such spaces is a class of 467 spaces each having 729 points and a constant line size of 8. These are a complete classification of point imprimitive, line transitive, linear spaces with this number of points and line size. The classification was completed in [7]. In fact, these are the only linear spaces achieving equality in the Delandtsheer-Doyen bound discussed in Section 4.4.

No more point imprimitive, line transitive, linear spaces have been identified in the last 9 years and a total of only 469 such spaces, excluding projective planes, have been found. This leads us to question whether such spaces are as rare as this suggests, or whether, with more research, many more such spaces might be identified.

## Conclusion

We have introduced here the concept of a line transitive, point imprimitive, finite, linear space. We have also noted that very few such spaces, apart from projective planes, are known. At this stage we have not introduced any methods for testing whether a given space is line transitive or point imprimitive, nor any way of constructing such spaces. These methods are based on a related area, design theory, which we will introduce in the following chapter.

# Chapter 4

# Design Theory

## Introduction

We now introduce the concept of a design. The problem of finding line transitive, point imprimitive, linear spaces is a special case of finding point imprimitive, line transitive designs. The idea of a design parameter set will be introduced. These parameters give information about the possible existence of such a design. The problem then becomes trying to show the existence of a design with parameters satisfying known conditions.

The advantage of searching for designs is that, although we are looking for a particular type of design, there are some results arising from more general design theory which can be used to these ends. Two of these results, the Orbit Theorem and the Delandtsheer-Doyen Theorem, are presented here. How these results can be used as part of a search will be discussed. This will be put into practise in Chapter 5.

## 4.1 2-$(v, k, 1)$ Designs

A 2-$(v, k, 1)$ design may be thought of as a special type of linear space. In essence, it is a linear space where all the lines contain the same number of points.

**Definition 4.1.1.** Suppose $\mathcal{D} = (\mathcal{P}, \mathcal{L})$ where $\mathcal{P}$ is a finite set and $\mathcal{L}$ is a set of subsets of $\mathcal{P}$. We refer to elements of $\mathcal{P}$ as *points* and of $\mathcal{L}$ as *lines* or *blocks*. Then $\mathcal{D}$ is a 2-$(v, k, 1)$ *design* if there are $v$ points, each line is a $k$-subset of $\mathcal{P}$ and, given any 2 points in $\mathcal{P}$, there is exactly 1 line in $\mathcal{L}$ containing them both.

We see that the last condition is exactly the same as the linear space condition discussed in Chapter 3. From now, when $v$ and $k$ are not relevant, we will refer to 2-$(v, k, 1)$ designs as designs.

2-$(v, k, 1)$ designs are a specific class of 2-$(v, k, \lambda)$ designs. The definition of a 2-$(v, k, \lambda)$

design is the same, except that every 2 points lie on exactly $\lambda$ lines.

Given the similarity between designs and linear spaces, it is not surprising that automorphisms are defined in the same way.

**Definition 4.1.2.** Let $\mathcal{D} = (\mathcal{P}, \mathcal{L})$ be a design. An *automorphism of* $\mathcal{D}$ is a map $\pi \in \mathrm{Sym}(\mathcal{P})$ such that $L^{\pi} \in \mathcal{L}$ for all $L \in \mathcal{L}$. A group $G \leq \mathrm{Sym}(\mathcal{P})$ is an *Automorphism group* of $\mathcal{D}$ if each of its elements is an automorphism of $\mathcal{D}$. The group of all automorphisms of $\mathcal{D}$ is denoted $\mathrm{Aut}(\mathcal{D})$.

Here, as before, we define the action of $G$ on lines to the induced action on $k$-subsets of points. Given the definitions are the same, automorphism groups are the same whether a space is considered as a design or just as a linear space.

We noted earlier that if a linear space is line transitive, then its line set is the orbit of any line under the automorphism group. It follows, since automorphisms are one to one, that the size of the lines is constant in a line transitive, linear space.

**Lemma 4.1.3.** *Let* $(\mathcal{P}, \mathcal{L})$ *be a linear space and* $G$ *be a line transitive automorphism group. Then the line size of the space is constant.*

**Proof.** Take some fixed $L \in \mathcal{L}$. Now take any $M \in \mathcal{L}$. Since $G$ is transitive on $\mathcal{L}$, there is some $g \in G$ such that $M = L^g$. It follows that $M = \{p^g | p \in L\}$ and hence, since $g$ is 1-1, $|M| = |L|$. This is true for any $M$, so the line size is constant. □

Since any linear space with a constant line size is a 2-$(v, k, 1)$ design for some $v$ and $k$, any line transitive linear space is a line transitive 2-$(v, k, 1)$ design. The line transitivity of the design follows since any automorphism group of the space considered as a linear space is also an automorphism group of the space considered as a design.

## 4.2 Parameter Sets

We are interested in identifying line transitive, point imprimitive linear spaces. Since we know all such spaces are designs, it is sufficient to identify $G$-point imprimitive, line transitive 2-$(v, k, 1)$ designs.

The first step in this is to identify $v$, $k$ and $G$ such that such a design might exist. We call $v$, $k$ and $G$ a *parameter set*. The need to identify a parameter set seems hardly worth mentioning. However, the approach to identifying parameter sets is very different to identifying designs.

To identify potential parameter sets, it is a matter of excluding those sets for which no design may exist. This is achieved using numerical and combinatorial arguments on $v$ and $k$ and using group theoretical arguments on $G$. Using these arguments it is possible to narrow the parameter sets to be considered down to a relatively small number. It is important to note however that this gives no guarantees that a design will exist with these parameters. It then becomes a matter of trying to construct a point imprimitive, line transitive design with the given parameters.

It should be noted that the group $G$ is not usually identified as a permutation group on $\mathcal{P}$. The possible permutation representations can be identified using group theoretic arguments on the information known about $G$. It is often the case that there may be more than one possible permutation representation for a given $G$, or several groups $G$ with specific detailed properties. Since there is no way of deciding if any one representation or group is more correct than another, each of these representations and groups needs to be searched separately.

The work of identifying parameter sets will not be considered here. In remainder of this dissertation we will assume such a set has already been identified.

## 4.3 The Orbit Theorem

We assume that we are working with a parameter set $v$, $k$, $G$, where $G$ has a prescribed permutation representation on some set $\mathcal{P}$ where $|\mathcal{P}| = v$. Since we require the final design to be point imprimitive, we must already have $G$ imprimitive on $\mathcal{P}$.

We wish to find $G$-line transitive 2-$(v, k, 1)$ designs $\mathcal{D} = (\mathcal{P}, \mathcal{L})$. The problem then is to identify one or more line sets $\mathcal{L}$ which result in a design. One brute force approach to solving this would be to search all possible sets of $k$-subsets of $\mathcal{P}$ and check which ones satisfy the design requirements and are orbits under the $G$ action. In almost every case this is computationally impossible.

We noted earlier that, since we require $G$ to be transitive on lines, $\mathcal{L}$ is the $G$-orbit of any line $L \in \mathcal{L}$. To identify a line set it is sufficient then to identify a single line $L$. Our problem reduces then to searching possible $k$-subsets of $\mathcal{P}$. It is still computationally impossible to check each set by first constructing its orbit under $G$ and then checking the design requirements. This is where our first design theory result becomes useful. The Orbit Theorem gives necessary and sufficient conditions on a $k$-subset for its orbit to be a design. In fact, this result was originally written for more general $t$-$(v, k, \lambda)$ designs, spaces where every $t$-subset of points lies on exactly $\lambda$ blocks, but we only consider it in the context useful to us here.

The ideas behind this theorem come from a construction used by Altop, outlined in [5]. This was used to construct 4-designs (similar to 2-designs, except every 4 points lie on exactly $\lambda$ lines), which were always 5-designs. The theorem as we present it here, along with the Delandtsheer-Doyen Theorem (Theorem 4.4.2) formed the basis for the search performed in [8].

**Theorem 4.3.1 (Orbit Theorem).** *Let $v$, $k \in \mathbb{N}$ and let $\mathcal{P}$ be a set with $|\mathcal{P}| = v$. Let $G$ be a permutation group on $\mathcal{P}$ and $O_1, O_2, \ldots, O_m$ be the orbits of $G$ on the 2-subsets of $\mathcal{P}$. Let $L$ be a $k$-subset of $\mathcal{P}$. Define $\hat{b} = \frac{v(v-1)}{k(k-1)}$ and for each $i = 1, 2, \cdots, m$, define*

*$b_i = |\{Y \in O_i : \ Y \subset L\}|$. Then there exists $\lambda \in \mathbb{N}$ such that $(\mathcal{P}, L^G)$ is a 2-$(v, k, \lambda)$ design if and only if the numbers $\frac{b_i}{|O_i|}$ are equal for all $i = 1, 2, \cdots, m$. Furthermore, if these numbers are all equal, then $b_i = \frac{|O_i|}{\hat{b}}$ for each $i$.*

The property that $\frac{b_i}{|O_i|}$ is independant of $i$ is called the *design property* and subsets satisfying it are called *starter blocks.* Note that the quantity $\hat{b}$ is not, in general, an integer. However, if $(\mathcal{P}, L^G)$ is a 2-$(v, k, 1)$ design, and hence a linear space, then $\hat{b}$ is, by Lemma 3.1.4, the number of lines in the space and hence an integer.

**Proof.** Suppose $\hat{L}$ is a $k$-subset of $\mathcal{P}$. Define $\mathcal{L} := \hat{L}^G = \{\hat{L}^g | \ g \in G\}$. It is clear that for any $L \in \mathcal{L}$, $|L| = k$ and $b_i = |\{Y \in O_i : \ Y \subset L\}|$, since G is transitive on $\mathcal{L}$ and $O_i$ is G-invariant. We wish to determine the number of blocks that contain a given 2-subset of $\mathcal{P}$. Take some $i$ with $1 \leq i \leq m$, then take $T_1, \ T_2 \in O_i$. Notice that since $O_i$ is G-transitive, there exists $g \in G$ such that $T_1^g = T_2$. It follows that $\{L \in \mathcal{L}| \ T_1 \subset L\}^g = \{L \in \mathcal{L}| \ T_2 \subset L\}$ and hence $|\{L \in \mathcal{L}| \ T_1 \subset L\}| = |\{L \in \mathcal{L}| \ T_2 \subset L\}|$. Thus the number of blocks containing a given 2-subset of $\mathcal{P}$ is constant within $O_i$. Define $\lambda_i := |\{L \in \mathcal{L}| \ T_1 \subset L\}|$ to be this number. Now examine the set of pairs $\{(L, T)| \ L \in \mathcal{L}, \ T \in O_i, \ T \subset L\}$. If we first select a block, then a 2-subset from $O_i$ contained in that block, we see there are $|\mathcal{L}| \cdot b_i$ such pairs. If we first select a 2-subset from $O_i$, then a block containing that subset, we see there are $|O_i| \cdot \lambda_i$ pairs. Thus

$$\lambda_i = |\mathcal{L}| \cdot \frac{b_i}{|O_i|}.$$

If $\hat{L}$ satisfied the design property, $\frac{b_i}{|O_i|}$ is independant of $i$ and hence $\lambda_i$ is independent of $i$. Thus we have a constant, call it $\lambda$, such that any 2-subset of $\mathcal{P}$ is contained in exactly $\lambda$ elements of $\mathcal{L}$. It follows that $(\mathcal{P}, \mathcal{L})$ is a 2-$(v, k, \lambda)$ design. Then, using the same argument used in Lemma 3.1.4, but noting that any 2 points are contained in $\lambda$ lines, we see that $\lambda\hat{b} = |\mathcal{L}|$. Hence $\frac{b_i}{|O_i|} = \frac{\lambda}{\lambda\hat{b}}$, and thus $b_i = \frac{|O_i|}{\hat{b}}$.

Conversely, suppose $(\mathcal{P}, \mathcal{L})$ is a 2-$(v, k, \lambda)$ design. Take $\hat{L}$ as above to be any block in $\mathcal{L}$. Clearly then $\lambda_i = \lambda$ for each $i$ and hence each $\frac{b_i}{|O_i|} = \frac{\lambda}{|\mathcal{L}|}$, which is independant of $i$. Notice by the above argument, this common value is $\frac{1}{\hat{b}}$. Thus any block satisfies the design property.

□

**Remark 4.3.2.** We wish to use this theorem to construct 2-$(v, k, 1)$ designs. Given a $k$-subset $L$ satisfying the design property, we must ensure that $\lambda = 1$ in the equations above. This gives us $|L^G| \cdot \frac{b_i}{|O_i|} = 1$. Since $b_i = \frac{|O_i|}{\hat{b}}$, this gives $|L^G| = \hat{b}$. This makes sense, since for a 2-$(v, k, 1)$ design, $\hat{b} = b$, the number of lines, as shown in Lemma 3.1.4. By the Orbit Stabiliser Theorem, this gives $|G_L| = \frac{|G|}{|L^G|} = \frac{|G|}{\hat{b}}$. This gives us a second condition to check on $k$-subsets. The design property combined with the extra condition that $|G_L| = \frac{|G|}{\hat{b}}$ is called the *orbit condition.*

Given $v$, $k$ and $G$ and a $k$-subset of $\mathcal{P}$, $L$, the Orbit Theorem is used in the following way to check the orbit condition on $L$. We first calculate $\hat{b} = \frac{v(v-1)}{k(k-1)}$. We use this to check that $|G_L| = \frac{|G|}{\hat{b}}$. If not, then $L$ cannot generate a 2-$(v,k,1)$ design. Next we construct $O_1, O_2, \cdots, O_m$, the orbits of $G$ on the 2-subsets of $\mathcal{P}$. Finally, for each $i = 1, 2, \cdots, m$, we calculate $b_i = \frac{|O_i|}{b}$ and check that there are exactly $b_i$ 2-subsets of $\mathcal{P}$ from $O_i$ contained in $L$ for each $i = 1, 2, \cdots, m$.

Notice that if we are checking many $k$-subsets, we only need to calculate $\hat{b}$, $b_1, b_2, \cdots, b_m$ and $O_1, O_2, \cdots, O_m$ once. Checking $|G_L|$ and the number of 2-subsets per $k$-subset per orbit are computationally simple, so this test can be performed very quickly.

## 4.4 Inner Pairs

Although the Orbit Theorem makes the check on $k$-subsets simple for a given $k$-subset, there are almost always more $k$-subsets than can reasonably be checked. As we will see in Chapter 5 there may be in excess of $10^{20}$ such subsets. This number cannot be searched given current computing power.

We need some further conditions to narrow the $k$-subsets we need to consider down further. The first condition arises from the work of Delandtsheer and Doyen in [3]. It is based on conditions all lines of a point imprimitive design must satisfy.

The Delandtsheer-Doyen Theorem primarily gives numerical conditions on $v$ and $k$ based on a non trivial G-invariant partition on the points. This numerical condition is used by researchers in identifying parameter sets.

As we see in the proof, one of the numbers involved in this condition has a significance to the design that we can use in our search. In order to see this, we introduce the concept of an inner pair.

**Definition 4.4.1.** Let $\mathcal{C} = \{C_1, C_2, \cdots, C_d\}$ be a non trivial partition on the points of $\mathcal{D} = (\mathcal{P}, \mathcal{L})$. Then an *inner pair* on a line $L \in \mathcal{L}$ is a 2-subset of points in $L$ contained within the same class $C_i$. Any pair of points that is not an inner pair is called an *outer pair*

The following theorem and its proof were the main result of [3].

**Theorem 4.4.2 (Delandtsheer-Doyen).** *Let $\mathcal{D} = (\mathcal{P}, \mathcal{L})$ be a 2-$(v,k,1)$ design and $G$ be an automorphism group of $\mathcal{D}$. Suppose there is a non-trivial $G$-invariant partition $\mathcal{C}$ of the point set $\mathcal{P}$ with $d$ classes, each of size $c$. Then there exist $x, y \in \mathbb{Z}^+$ such that*

$$c = \frac{\binom{k}{2} - x}{y} \text{ and } d = \frac{\binom{k}{2} - y}{x}. \tag{4.4.1}$$

*As a result, since $v = cd$,*

$$v \leq (\binom{k}{2} - 1)^2 \tag{4.4.2}$$

**Proof.** Since $G$ is line transitive, the number of inner pairs per line, and hence the number of outer pairs per line, is constant. Let $x, n$ denote the number of inner and outer pairs per line respectively. Then

$$x + n = \binom{k}{2} \tag{4.4.3}$$

and $x, n \geq 1$ since any pair is contained in at least one line. Consider the pairs $(\{p_1, p_2\}, L)$ of an outer pair $p_1, p_2 \in \mathcal{P}$ and a line $L \in \mathcal{L}$ with $p_1, p_2 \in L$. By counting the points first then the line, and the line first then the points, we see that, if $b = |\mathcal{L}|$,

$$bn = \binom{d}{2} c^2. \tag{4.4.4}$$

We know from Lemma 3.1.4, since $v = cd$, that

$$b = \frac{cd(cd-1)}{k(k-1)}. \tag{4.4.5}$$

From Equations 4.4.3, 4.4.4 and 4.4.5 we see

$$\begin{aligned}
ncd(cd-1) &= k(k-1)\binom{d}{2}c^2 \\
nd(cd-1) &= 2\binom{k}{2}\frac{d(d-1)}{2}c \\
n(cd-1) &= \binom{k}{2}(d-1)c \\
n(cd-1) &= (x+n)cd - \binom{k}{2}c \\
-n &= xcd - \binom{k}{2}c
\end{aligned}$$

and hence

$$\binom{k}{2} - dx = \frac{n}{c}. \tag{4.4.6}$$

Set $y = \binom{k}{2} - dx$. Then clearly $y \in \mathbb{Z}^+$. Since $n = \binom{k}{2} - x$, clearly $c = \frac{\binom{k}{2}-x}{y}$. Also, since $y = \binom{k}{2} - dx$ we have $dx = \binom{k}{2} - y$ and hence $d = \frac{\binom{k}{2}-y}{x}$. □

We see then from the proof that $x$ in Equation 4.4.2 is the number of inner pairs per line, with respect to the partition $\mathcal{C}$.

The number of inner pairs tells us something about how many points there are from a given line per class. Since there is nothing to distinguish the classes in this information, it can at most tell us how many classes contain a given number of points, not which class. We must also note that this information is not necessarily unique. For example, if $x = 1$ then there is a unique inner pair. The only way this can happen is if, for any given line, there is exactly one class containing 2 points and all other classes contain either 1 or 0 points from the line. However, if $x = 3$ it may be that, for a given line, 3 classes contain 2 points each and all other classes contain 0 or 1, or it may be that 1 class contains 3 points and all others contain 0 or 1 point. Notice 3 points in the same class gives $\binom{3}{2} = 3$ possible 2-subsets and hence 3 inner pairs.

We formalise this information as an intercept vector for the lines. This vector can be thought of as a list $[d_0, d_1, \cdots, d_n]$ where, for any line and each $i$, $d_i$ is the number of classes containing $i$ points from that line. Since $G$ is transitive on the line set and preserves the partition, this vector will be constant across all lines. As with the permutation representation of $G$, this vector may not be uniquely determined by the parameters $(c, d, x, y, G)$, so we may need to search more than once, using a different intercept vector each time.

## 4.5 Symmetry

There is no more information gained by identifying a number of isomorphic designs than by identifying one of those designs. Once a design has been identified, we could if we wished construct a number of isomorphic designs from it as a separate exercise to the search. In order to reduce the amount of work that must be done during the search, we only try to search for k-subsets that will yield non-isomorphic designs.

We must be careful however to only consider isomorphisms that preserve the automorphism group $G$ and its line transitivity and point imprimitivity. It may be that a design with $G$ as an automorphism group has an isomorphic design that does not admit $G$. In this case, the isomorphic design does not conform to our search requirements and should not be considered.

We may define an equivalence relation on $k$-subsets of $\mathcal{P}$ that satisfy the design property of Theorem 4.3.1 based on isomorphism of designs. We will say two such sets $L_1$ and $L_2$ are equivalent if there exists some $\pi \in \mathrm{Sym}(\mathcal{P})$ such that $(L_1^G)^\pi = L_2^G$. If $L_1$ and $L_2$ are equivalent and both starter blocks, then it follows, by the definition of an isomorphism, that the designs $(\mathcal{P}, L_1^G)$ and $(\mathcal{P}, L_2^G)$ are isomorphic.

Note that given two $k$-subsets $L_1$ and $L_2$, it is very difficult to decide if there is some element of $\mathrm{Sym}(\mathcal{P})$ mapping $L_1^G$ to $L_2^G$. We instead work with subgroups of $\mathrm{Sym}(\mathcal{P})$.

Clearly then, we would like to check as few $k$-subsets per equivalence class as possible. We would like to identify some property that only one, or at least very few, $k$-subset per equivalence class has.

We first note that if two $k$-subsets are in the same $G$-orbit, then the line sets generated by them will be equal and hence the $k$-subsets are equivalent. We thus only need to check one $k$-subset per $G$-orbit. At the basic level, we could note that $G$ is transitive on points. We could pick some arbitrary point $\alpha$. Given any $k$-subset $L$, we could take any point $\beta \in L$ and there would be some $g \in G$ such that $\alpha = \beta^g$. Hence $\alpha \in L^g$. Thus we only need to consider $k$-subsets containing $\alpha$, since every equivalence class contains at least one such $k$-subset. This is not using as much information as we have however. Since, if $(\mathcal{P}, \mathcal{L})$ is not a projective plane and $G$ is line transitive, $G$ will never be regular on $\mathcal{P}$, we may be able to fix $\beta$ and move the remaining points around in some way. Ideally we would like to fix as many points as possible, to narrow down as much as possible the number of $k$-subsets to consider.

This argument has so far only considered $k$-subsets generating equal line sets, which is much more restrictive than than 'isomorphic' line sets. We would like to use a similar argument with a group of isomorphisms of the designs. The following lemma gives a method of identifying such a group.

**Lemma 4.5.1.** *Let $G$ be point imprimitive and block transitive on a 2-$(v,k,1)$ design $\mathcal{D} = (\mathcal{P}, \mathcal{L})$. Let $H$ be such that $G \trianglelefteq H \leq \mathrm{Sym}(\mathcal{P})$. Then for any $h \in H$, the image of $\mathcal{D}$ under $h$, $\mathcal{D}^h = (\mathcal{P}, \mathcal{L}^h)$ is a 2-$(v,k,1)$ design isomorphic to $\mathcal{D}$. Furthermore, $G$ is an automorphism group of $\mathcal{D}^h$, and is point imprimitive and line transitive on $\mathcal{D}^h$.*

**Proof.** Let $\mathcal{M} = \mathcal{L}^h$. Then $\mathcal{M}$ is a set of $k$-subsets of $\mathcal{P}$, since $h$ is 1-1. Also, since $h$ is 1-1 on points, it must be 1-1 on 2-subsets. Thus any two points $p_1, p_2 \in \mathcal{P}$ are the image of a unique pair $q_1, q_2$. There is a unique line $L \in \mathcal{L}$ containing $q_1$ and $q_2$ so it follows that $L^h$ is the unique line in $\mathcal{M}$ containing $p_1$ and $p_2$. Hence $\mathcal{D}^h$ is a 2-$(v,k,1)$ design. Since $h \in \mathrm{Sym}(\mathcal{P})$ it is clear that $\mathcal{D}^h \cong \mathcal{D}$.

We now need to show $G \leq \mathrm{Aut}\mathcal{D}^{\mathrm{h}}$ and $G$ is point imprimitive and block transitive. Take $M \in \mathcal{M}$ and $g \in G$. Then $M = L^h$ for some $L \in \mathcal{L}$. Hence $M^g = L^{(hg)}$. Since $G \trianglelefteq H$, we know $hG = Gh$ so $hg = g_1 h$ for some $g_1 \in G$. Thus $M^g = (L^{g_1})^h \in \mathcal{L}^h = \mathcal{M}$. So $G$ maps lines to lines. Hence $G$ is an automorphism group on $\mathcal{D}^h$. Line transitivity follows from the same argument. Take $M_1, M_2 \in \mathcal{M}$. Then $M_1 = L_1^h, M_2 = L_2^h$ for some $L_1, L_2 \in \mathcal{L}$. Since $G$ is transitive on $\mathcal{L}$, there is some $g \in G$ such that $L_1^g = L_2$ and since $gh = hg_1$ for some $g_1 \in G$ it follows $M_1^{g_1} = M_2$. Finally, $h$ preserves the point set, so $G$ is still point imprimitive. □

So if we can find some subgroup of $\mathrm{Sym}(\mathcal{P})$ containing $G$ as a normal subgroup, we have a group of isomorphisms that preserve the transitivity and imprimitivity of $G$. We then need only search one $k$-subset from each orbit of this new, larger group.

We refer to this process of narrowing down the subsets to search as a symmetry argument, since we are looking at orbits under a group action.

# Conclusion

There are two steps in searching for point imprimitive line transitive linear spaces. We have shown all such spaces must be 2-$(v,k,1)$ designs, so the first step is identifying potential parameters $v$ and $k$ and a group $G$. The second set is trying to identify a design with these parameters.

We have seen the search for a design becomes a search for $k$-subsets satisfying the conditions given by the Orbit Theorem and the Delandtsheer-Doyen Theorem. Such a search would normally be impossible, due to the large number of possible subsets. By identifying an overgroup of $G$ that normalises $G$, we can apply a symmetry argument to reduce the number of $k$-subsets to examine by a factor of the order of this group.

# Chapter 5

# Examples of Two Searches

## Introduction

In this chapter we will introduce two parameter sets identified by other researchers. We will outline how the methods given in Chapter 4 are applied to these parameter sets. The information gained here will be applied in Chapter 6 to construct a search algorithm. Although neither of the sets yielded designs, the methods used here can be applied to a wider class of parameters.

## 5.1 The Parameter Sets

The first parameter set we examined arose from [10, Theorem 6]. This was one of a very small number of parameter sets satisfying an extra group theoretic condition, and having parameter $x$ from Theorem 4.4.2 small. It has $v = 451 = 41 \cdot 11$ and $k = 10$. By this we mean the point set has 451 points and each line contains 10 points. As is often the case, the group $G$ acting on the points is only very generally described. However, we are given that there is a $G$ invariant partition with 11 classes of 41 points each. Further, the subgroup $K$ of $G$ which fixes, set-wise, every class $C$ of this partition is isomorphic to $D_{82}$, the dihedral group on 41 points. The dihedral group can be thought of as the group of symmetries - rotations and reflections - of a 41-gon. A little more information about $G$ is given, which will be discussed in the next section. We will call this parameter set 1.

The second parameter set is very similar to set 1. After the search on set 1 failed, the logical step is to try to identify what part of the parameter information may have prevented it from yielding a design. As we will see in section 5.5, the fact that $D_{82}$ acts irregularly on the class $C$ restricts the possible starter blocks greatly. Recall the definition of a regular action from Section 2.3. If we change $G$ so that $K$ acts regularly on $C$, it may increase the possibilities for valid starter blocks. Recall a starter block, as described in Section 4.3 is a 10-subset of the point-set which satisfies the design property given in

Theorem 4.3.1.

Parameter set 2 thus begins with similar assumptions to set 1. We have $v = 451$ and $k = 10$. We suppose there is a $G$-invariant partition with 11 classes of size 41, such that the subgroup $K$ fixing every class $C$ of this partition set-wise is isomorphic to $\mathbb{Z}_{41}$, the cyclic group with 41 elements. This group acts regularly on each class $C$. Once this relaxation of the conditions of parameter set 1 had been made, it was noticed that parameter set 2 had been identified as a possible parameter set by the work of a group at UWA[1]. Since it had been identified by two separate processes, parameter set 2 seemed an interesting set to search.

We now have the basic information to start the search process for two searches. As we saw in the previous chapter, there are several steps to perform before a computer based search for designs can be done.

## 5.2 The Permutation Group $G$

The first important step in any design search is identifying explicit possibilities for the group $G$ and its permutation representations on $\mathcal{P}$. Without these representations, the Orbit Theorem criteria cannot be applied, since it relies on $G$-orbits. The work of identifying these possible representations was performed by Niemeyer and Praeger and an outline will be given here.

The first step performed by Niemeyer and Praeger [9] on both parameter sets was showing there there is a second $G$-invariant partition, having 41 classes of size 11. The existence of one of these two partitions, under a line transitive group, implies the existence of both, so we can always use the information given by either. The partitions are orthogonal, in the sense that the intersection of two classes from the different partitions is a unique point. We can thus think of these partitions as forming a grid structure over the points, having 11 columns and 41 rows. The first partition is the columns of the grid, the second is the rows of the grid. We can see this in Figure 5.2.1.

We require some labelling of the rows and columns to describe actions on them. Label the rows $0, 1, \cdots, 40$ and the columns $0, 1, \cdots, 10$. We will refer to the point in row $e$, column $f$ as $(e, f)$. At this stage there is nothing to distinguish different rows or different columns, so we may make this labelling arbitrarily. This is an identification of the points with $\mathbb{Z}_{41} \times \mathbb{Z}_{11}$. Performing this identification is sensible, since we will see later that $\mathbb{Z}_{41} \times \mathbb{Z}_{11}$ acts regularly on the points. We will often require a total ordering on the points. We will compare points by column first, then row. Define $(e, f) \leq (g, h)$ if $f < h$ or $f = h$ and $e \leq g$.

Having this grid structure makes describing the $G$ action easier, since we can describe the action of a group element by the way it permutes the rows and the columns. For both parameter sets it is clear then that $G \leq \mathrm{Sym}(41) \times \mathrm{Sym}(11)$ since $G$ preserves the rows and columns. In fact, it was shown [9] that for both parameter sets, $G \leq H = \mathrm{AGL}(1, 41) \times \mathrm{AGL}(1, 11)$.

[1]This group was Betten, Delandtsheer, Niemeyer and Praeger. No papers have yet resulted from this work.

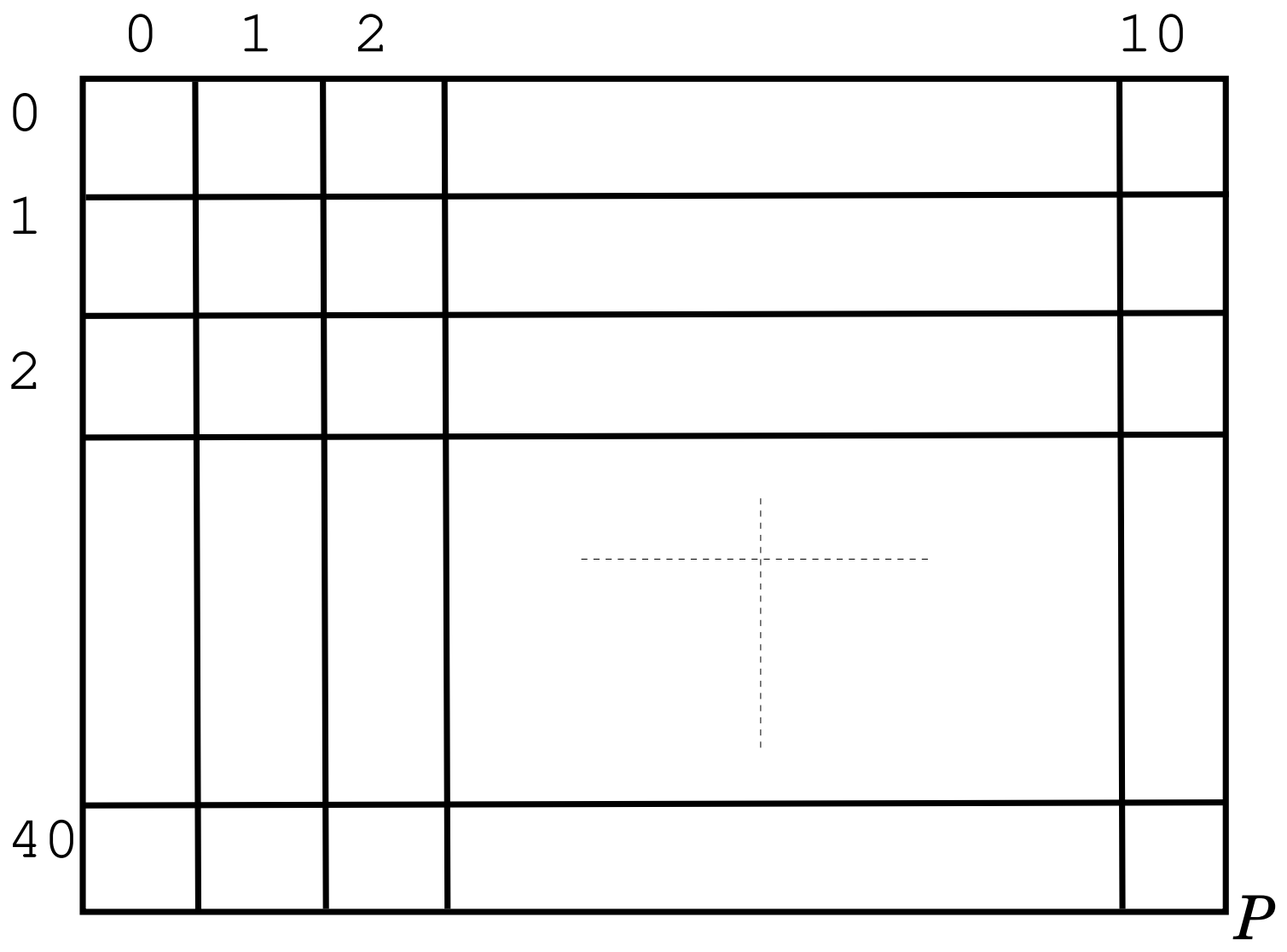

Figure 5.2.1: The grid structure on the points.

The Affine General Linear group on $\mathbb{Z}_p$, $\mathrm{AGL}(1, p)$ consists of all affine maps $\mathbb{Z}_p \to \mathbb{Z}_p$. These are maps of the form $\gamma \mapsto a\gamma + b$ for some $a, b \in \mathbb{Z}_p, a \neq 0$. Notice then that $\mathrm{AGL}(1, p) \cong \mathbb{Z}_p \cdot \mathbb{Z}_{p-1}$. In this case we have $\mathrm{AGL}(1, 41)$ acting on the rows and $\mathrm{AGL}(1, 11)$ acting on the columns of our grid structure. If $a$ and $b$ are primitive elements of $\mathbb{Z}_{41}$ and $\mathbb{Z}_{11}$ respectively, we can construct $H$ from the following generators:

- $\alpha$, addition by 1 to the row, so $\alpha : (e, f) \mapsto (e + 1, f)$.
- $\beta$, addition by 1 to the column, so $\beta : (e, f) \mapsto (e, f + 1)$.
- $a$, multiplication on rows, so $a : (e, f) \mapsto (ae, f)$.
- $b$, multiplication on columns, so $b : (e, f) \mapsto (e, bf)$.

In both parameter sets, $G$ contains, as a normal subgroup, the group of additions on the rows and columns, $N = \mathbb{Z}_{41} \times \mathbb{Z}_{11}$. We can see that $N = \langle \alpha, \beta \rangle$. Notice, as noted above, that $N$ acts regularly on $\mathcal{P}$, since $\mathbb{Z}_{41}$ acts regularly on the rows and $\mathbb{Z}_{11}$ on the columns. Hence identifying the points with $N$ is possible, as we saw in Lemma 2.3.4.

Niemeyer and Praeger showed that for parameter set 1, $G \cong (D_{82} \times \mathbb{Z}_{11}) \cdot \mathbb{Z}_5 \cong N \cdot (\mathbb{Z}_2 \times \mathbb{Z}_5)$ and for parameter set 2, $G \cong (\mathbb{Z}_{41} \times \mathbb{Z}_{11}) \cdot \mathbb{Z}_5 \cong N \cdot \mathbb{Z}_5$. Although this describes the overall structure of the group, it does not define how $G$ acts on the points.

In order to construct $D_{82}$ we need an involution on the rows. This will be a multiplicative element of order 2. Such an element is $a^{20}$, since $a$ has multiplicative order 40. Notice that $a^{20} \cong -1 \pmod{41}$.

In both cases, Niemeyer and Praeger showed that $\mathbb{Z}_5$ must act non-trivially on both that rows and the columns, and that its permutation representation would be some subgroup of $\langle a^8 \rangle \times \langle b^2 \rangle$. Notice $a^8$ and $b^2$ are multiplicative elements of order 5 of $\mathbb{Z}_{41}$ and $\mathbb{Z}_{11}$ respectively. There are six possible subgroups of this group of order 5. They are $\langle a \rangle, \langle ab \rangle, \langle ab^2 \rangle, \langle ab^3 \rangle, \langle ab^4 \rangle$ and $\langle b \rangle$. Notice that every of element $\langle a \rangle$ fixes every column and every element of $\langle b \rangle$ fixes every row. These groups are thus ruled out by the non-triviality condition. We are left with four possibilities for the action $\mathbb{Z}_5$, this leads to 4 possible, non-isomorphic, groups. We will denote the group constructed using $\langle a^8 b^{2i} \rangle$ by $G_i$ for $i = 1, 2, 3, 4$. Further, Niemeyer and Praeger showed that $G_i \ncong G_j$ for $i \neq j$. We will not give a proof of this non-isomorphism here, but it can be seen in Section 6.4 that the four different groups give different results during testing of our final algorithm. Isomorphic groups would have given the same results.

Let us consider how to construct $G_i$ for parameter set 1. Once again, let $a$ and $b$ be primitive elements of $\mathbb{Z}_{41}$ and $\mathbb{Z}_{11}$ respectively. Then $G_i$ has 4 generators:

- $\alpha$, addition by 1 to the row, so $\alpha : (e, f) \mapsto (e + 1, f)$.
- $\beta$, addition by 1 to the column, so $\beta : (e, f) \mapsto (e, f + 1)$.
- $a^{20}$, the involution on rows, so $a^{20} : (e, f) \mapsto (-e \pmod{41}, f)$.
- $a^8 b^{2i}$, the order 5 multiplication on rows and columns where $a^8 b^{2i} : (e, f) \mapsto (a^8 e, b^{2i} f)$.

Constructing $G_i$ for, parameter set 2 is very similar. The difference between the two parameter sets is the involution on rows, $a^{20}$. Hence the generators of $G_i$ for parameter set 2 are $\alpha$, $\beta$, and $a^8 b^{2i}$.

We denote by $N_{\mathrm{Sym}(\mathcal{P})}(G)$ the *normaliser of* $G$ *in Sym*$(\mathcal{P})$. This is the largest subgroup of $\mathrm{Sym}(\mathcal{P})$ that contains $G$ as a normal subgroup, and if $G \trianglelefteq M \leq \mathrm{Sym}(\mathcal{P})$, then $M \leq N_{\mathrm{Sym}(\mathcal{P})}(G)$.

**Lemma 5.2.1.** *The normaliser of G in Sym*$(\mathcal{P})$ *is* $N_{\mathrm{Sym}(\mathcal{P})}(G) = H$.

**Proof.** Since $N \trianglelefteq G$ and $N$ is the unique normal subgroup of order 451, any element of $\mathrm{Sym}(\mathcal{P})$ which normalises $G$ must normalise $N$. It follows that $N_{\mathrm{Sym}(\mathcal{P})}(G) \leq N_{\mathrm{Sym}(\mathcal{P})}(N)$. By [4, Exercise 2.5.6, p45] we know that the normaliser of a subgroup $M$ of $\mathrm{Sym}(\mathcal{P})$ that is regular on $\mathcal{P}$ is $M \cdot \mathrm{Aut(M)}$. We claim that $\mathrm{Aut(N)} = \mathbb{Z}_{40} \times \mathbb{Z}_{10}$.

We know that $\mathrm{Aut}(\mathbb{Z}_{41}) = \mathbb{Z}_{40}$ and $\mathrm{Aut}(\mathbb{Z}_{11}) = \mathbb{Z}_{10}$, both acting by multiplication. Given $\phi$, an automorphism of $\mathbb{Z}_{41}$ and $\psi$, an automorphism of $\mathbb{Z}_{11}$, we can construct an automorphism of $\mathbb{Z}_{41} \times \mathbb{Z}_{11}$ by $(e, f) \mapsto (\phi(e), \psi(f))$. Hence $\mathbb{Z}_{40} \times \mathbb{Z}_{10} \subseteq \mathrm{Aut(N)}$. Now suppose $\varphi \in \mathrm{Aut(N)}$. Since $\varphi$ is a group automorphism, it must map any normal subgroup to another, isomorphic normal subgroup. It follows, since $\gcd(41,11)=1$, that $\varphi(\mathbb{Z}_{41})=\mathbb{Z}_{41}$ and $\varphi(\mathbb{Z}_{11})=\mathbb{Z}_{11}$. Hence $\varphi|_{\mathbb{Z}_{41}}\in \mathrm{Aut}(\mathbb{Z}_{41})$ and $\varphi|_{\mathbb{Z}_{11}}\in \mathrm{Aut}(\mathbb{Z}_{11})$. It remains to show that $\varphi=(\varphi|_{\mathbb{Z}_{41}},\varphi|_{\mathbb{Z}_{11}})$. Notice $\varphi(e,f)=\varphi((e,1)(1,f))=(\varphi|_{\mathbb{Z}_{41}}(e),1)(1,\varphi|_{\mathbb{Z}_{11}}(f))=(\varphi|_{\mathbb{Z}_{41}}(e),\varphi|_{\mathbb{Z}_{11}}(f))$ as required. Hence $\varphi\in\mathbb{Z}_{40}\times\mathbb{Z}_{10}$, so $\mathrm{Aut(N)}\subseteq\mathbb{Z}_{40}\times\mathbb{Z}_{10}$.

We see then that $N_{\mathrm{Sym}(\mathcal{P})}(G)\leq N\cdot(\mathbb{Z}_{40}\times\mathbb{Z}_{10})=H$. We will show equality by proving $G\trianglelefteq H$, since then it follows that $H\leq N_{\mathrm{Sym}(\mathcal{P})}(G)$. It is known that if $N\trianglelefteq G$ and $G\leq H$ then $G\trianglelefteq H$ if $\frac{G}{N}\trianglelefteq\frac{H}{N}$. Clearly $N\trianglelefteq G$, and notice $\frac{H}{N}\cong\mathbb{Z}_{40}\times\mathbb{Z}_{10}$, which is abelian. Any subgroup of an abelian group is normal, so it follows then that $\frac{G}{N}\trianglelefteq\frac{H}{N}$ and hence $G\trianglelefteq H$. □

We have now identified, for both parameter sets, all possible groups $G$ which arise as permutation groups on $\mathcal{P}$. Before applying the Orbit Theorem (Theorem 4.3.1), we will use the Delandtsheer-Doyen Theorem (Theorem 4.4.2) to examine the overall structure of a line in any possible design. We will see in Section 5.5 that this structural information will assist in applying the Orbit Theorem.

## 5.3 Applying the Delandtsheer-Doyen Theorem

Since the Delandtsheer-Doyen Theorem (Theorem 4.4.2) only uses information about the $G$-invariant partitions, there is no difference between its application to parameter set 1 and parameter set 2. We will give here the argument used on both parameter sets.

On these parameter sets, the Delandtsheer-Doyen Theorem will be particularly useful, since there are two different partitions to apply it to. We will thus be given inner pair information for both the partitions. We will see in Section 5.4 that this information can be put together to greatly restrict the possible $k$-subsets we need to examine. We will call the partition with 11 classes of size 41 (the columns) partition 1. Using the notation of Theorem 4.4.2, we have $d_1=11, c_1=41$. The second partition, the rows, we will call partition 2 and this gives us $d_2=41, c_2=11$.

The application of this theorem will give us information about the number of points from a given line in the same column or row. We will use the notation $n$-row to denote a row containing $n$ points from a given line and $m$-column to denote a column containing $m$ points from a given line.

We wish to calculate the values for $x$ and $y$ in Equation 4.4.2 for both partitions. We call these $x_1, y_1$ and $x_2, y_2$ for the first and second partitions respectively. Since $c_1=d_2, c_2=d_1$ it is clear that $x_2=y_1, y_2=x_1$. We are only interested in $x$ values, so we will solve equation 4.4.2 for partition 1 to get $x_1$ and $y_1$ and hence $x_2$.

This gives us $41=\frac{\binom{10}{2}-x_1}{x_2}$ and $11=\frac{\binom{10}{2}-x_2}{x_1}$, solving to $x_1=4$, $x_2=1$. We thus know there are 4 inner pairs on the first partition, the columns, and 1 inner pair on the second partition, the rows. This gives us possible intercept vectors of $[4,5,1,1]$ or $[3,4,4,0]$ for the columns and $[32,8,1,0]$ for the rows. Recalling the definition of these vectors from Section 4.4 we can interpret these. There is only one possibility for the rows. Given any line, there will be thirty two 0-rows, containing no points from the line, eight 1-rows containing 1 point from the line, and the remaining row will contain 2 points from the line. There are two possibilities we must consider independently for the columns. First it may be that given any line, there will be four 0-columns containing no points, five 1-columns, one 2-column and one 3-column. The other possibility is that given any line, there will be three 0-columns, four 1-columns and four 2-columns. It is important to recall that all lines in the design will have the same intercept vector. We will thus search separately for designs having the first intercept vector for columns and designs having the second intercept vector for columns.

It is important to consider how the intercept vectors on rows and columns will work together. For example, when we have a column containing 3 points from a line, how many points do the three rows containing these three points each contain? For each new point along these rows, how many points does its column contain? To answer these questions, we introduce the idea of a mask in Section 5.4.

## 5.4 Masks

A mask can be thought of as a general map of what any line in the design will look like. It is a picture of the line, with all points included, but without the rows or columns being labelled. From the mask it is clear how many rows or columns contain any given number of points from the line, and how these rows and columns interact. The search can then be run to find designs conforming to each different mask individually. At present these masks are constructed manually and any information they provide is given as input to the search. It is hoped that in the future a program will be developed to perform the mask construction automatically, but this is outside the scope of this dissertation.

Obviously there will be different masks for the two possible intercept vectors on columns. Once we have constructed all possible masks across all possible intercept vectors, there is no further need to differentiate between the various intercept vectors. Different intercept vectors are distinguished by the different masks they produce.

In this search the construction of the masks is made simpler since there is only one 2-row. Once we have decided how many points are in the columns containing the two points from the 2-row, all other points must be on 1-rows. This will completely determine the mask. The set of all possible masks for our parameter sets can be seen in Figure 5.4.1.

We consider first the intercept vector $[4, 5, 1, 1]$ on the columns. There are 4 possibilities for the 2-row. Firstly, its two points may each lie on 1-columns. This gives us mask 1A. Second, one of the points in the 2-row might lie of a 2-column and the other on a 1-column. This gives us mask 1B. Next, one point may be on the 3-column, the other on a 1-column. This is mask 1C. Finally, one of the points may be on the 3-column, the other of them a 2-column, which is mask 1D.

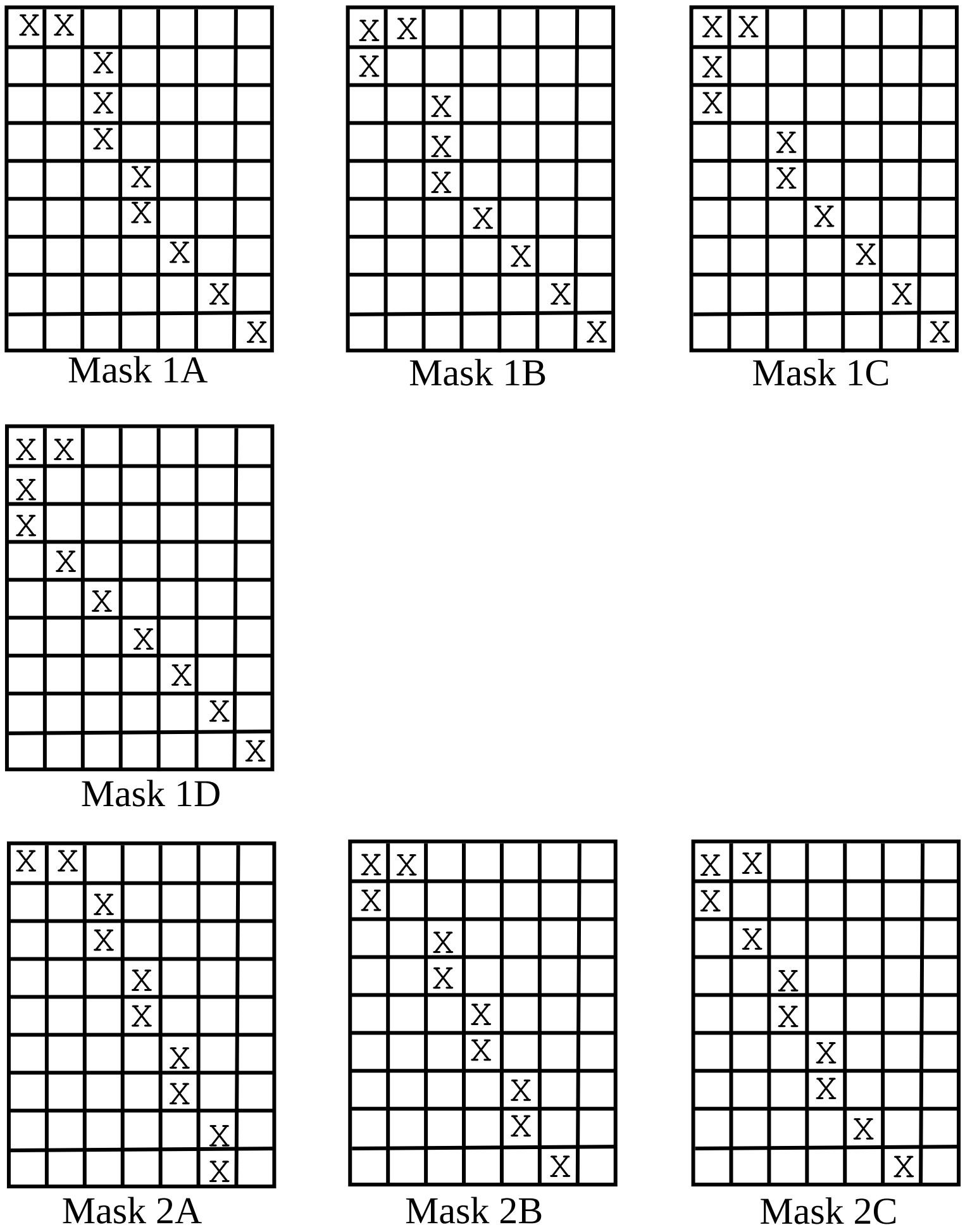

Figure 5.4.1: The possible masks for either parameter set

We now consider the second possible intercept vector, $[3, 4, 4, 0]$ on columns. Consider again the two points contained in the unique 2-row. It may be that both points are on 1-columns, which is mask 2A. Second, one point may be on a 1-column, the other on a 2-column, giving mask 2B. Finally, both points could lie on 2-columns, giving mask 2C.

We see then that there are 7 possible masks for this search. This narrows down the possibilities for starter blocks to be considered considerably. Next we will examine the Orbit Theorem's application and a symmetry argument to narrow down these possibilities further.

## 5.5 Using the Orbit Theorem

We will need to treat the two parameter sets separately when applying the Orbit Theorem (Theorem 4.3.1). This follows since the theorem involves the $G$-orbits of pairs of points and the groups involved are different. Here we will give a description of the $G$-orbits on 2-subsets of points and see how many 2-subsets a starter block must have from each orbit in each case. We will also see restrictions the theorem places on the starter block in parameter set 1.

We begin by applying of the theorem to parameter set 1. Recall that we must consider four groups $G$ separately. We consider the $G$-orbits on 2-subsets of points. We know that $G$ preserves the rows and columns, it follows that there will be 3 different types of orbits. There will be orbits where the points in the 2-subsets come from the same row but different columns, orbits where the points come from the same column but different rows and orbits where the points come from different rows and columns. Given the small numbers involved in this example, the orbits can be constructed computationally. The GAP code used to do this computation is found in Appendix A. From this computation we find that, for each of the four groups $G$, there are 25 orbits on 2-subsets. Of these, 5 orbits, those containing pairs of points from the same column or the same row, each contain $\frac{|G|}{2} = 2255$ pairs of points. The remaining 20 orbits each contain $|G| = 4510$ pairs. These are the 2-subsets from different rows and columns.

We are looking for linear spaces, so we require $\hat{b} = b$, using $\hat{b}$ as defined in the Orbit Theorem and $B$ the number of lines. Recall $\hat{b} = \frac{v(v-1)}{k(k-1)} = \frac{451 \cdot 450}{10 \cdot 9} = 451 \cdot 5 = \frac{|G|}{2}$, it follows from the Orbit Theorem (Theorem 4.3.1) that each starter block contains $\frac{|O_i|}{b} = \frac{2|O_i|}{|G|}$ 2-subsets from orbit $O_i$. Thus a starter block must contain exactly one 2-subset from each of the first 5, smaller orbits and two 2-subsets from the remaining 20 orbits. This gives us the orbit condition on 10-subsets to be checked during the search. Moreover, we note that $|G_L| = \frac{|G|}{b} = 2$. We use this fact to gain further information about the masks.

**Lemma 5.5.1.** *For parameter set 1, the intercept vector on columns is [5,2,4,0] and the mask is mask 2A of Figure 5.4.1.*

**Proof.** Let $L$ be a valid starter block for a 2-$(v, k, 1)$ design corresponding to parameter set 1. Note that a pair of points will uniquely determine a line. Thus if $L$ is fixed, the two points contained in its unique 2-row must also be fixed. Call this 2-subset $S$. There are two elements of $G$ that will fix $S$, by the Orbit Stabiliser Theorem (Theorem 2.3.3), since we know $|S^G| = \frac{|G|}{2}$. We can see that one of these elements must be the identity, 1. The other will be an involution from the normal subgroup $D_{82}$. We know every involution in $D_{82}$ fixes exactly one row and permutes the others in pairs. Let $j$ be the involution which fixes the row containing $S$. We thus require the element $j \in (D_{82} \times \mathbb{Z}_{11}) \cdot \mathbb{Z}_5$ to fix $L$ set-wise. We know it fixes the row containing $S$ and permutes all other rows, but preserves all columns. Hence $j$ must fix all 11 points in the row containing $S$. With the exception of the columns containing $S$ there can thus be no 1-columns. If there were, the point in this column would be mapped to a different point in the same column, which could not be in $L$, so $L$ would not be fixed. This tells us $L$ can have at most two 1-columns. Hence its intercept vector on the columns must be $[5, 2, 4, 0]$. Also, since these two 1-columns must make up the 2-row, we are restricted to only using the mask 2A of Figure 5.4.1. □

We have thus restricted the search on parameter set 1 to only one of its seven previously possible masks. This argument tells us more than just which mask we must use. We noted that $j \in (D_{82} \times \mathbb{Z}_{11}) \cdot \mathbb{Z}_5$ must fix the line $L$. For this element to fix $L$, the rows of the 2 points in any given 2-column are restricted. We need these rows to be one of the pairs that $j$ alternates. There are 20 possible pairs of rows and each of the four 2-columns must use a unique pair from these 20.

We now look at parameter set 2. Applying the Orbit Theorem here is much easier. From computational work in GAP, we find that, for each of the four possible groups $G_i$, there are 45 orbits on 2-subsets of points, each containing 2255 subsets.In this case, $b = 2255 = |G|$ and the Orbit Theorem implies that a starter block contains exactly $\frac{|O_i|}{b} = \frac{|G|}{|G|} = 1$ 2-subsets from each $G$-orbit $O_i$ on 2-subsets. Further, we require that for any starter block $L$, $|G_L| = 1$. As above, not that if $L$ is fixed then $S$, the 2 points it contains in its unique 2-row must also be fixed. However, since the size of each orbit on 2-subsets is $|G|$, it follows that $|G_S| = 1$ and hence $G_L| = 1$ as required. We need place no further restrictions on starter blocks to ensure that $\lambda = 1$.

We have seen here how the Orbit Theorem is applied to our two parameter sets. In each case, the information on the number of 2-subsets per orbit per line was easily obtained. In addition, a tight restriction was found in Lemma 5.5.1 on starter blocks for parameter set 1. This restriction means the search will be much easier and quicker to perform, but also means the parameter set is much less likely to produce a design. As we will see, this parameter set did not yield a design. The change between parameter sets 1 and 2 was designed to remove this restriction. The search on that parameter set will thus be longer, but may be more likely to yield a design.

## 5.6 Symmetry Arguments

So far in this chapter, we have built up necessary and sufficient conditions for a given 10-subset of $\mathcal{P}$ to generate a 2-$(451, 10, 1)$ design under our given group $G$. However, in both cases, the number of possible $k$-subsets we need to check is still far too large for a search to be performed in realistic time. As we noted in Section 4.5, we can greatly reduce the number of 10-subsets we need to consider by only considering one starter block per equivalence class. We also saw that equivalence classes can be generated by considering the orbits of $k$-subsets under some group $M$ where $G \trianglelefteq M$ and $M$ preserves the $G$-invariant partition. Here we will construct such a group for both searches and show how its orbit information may be used to restrict the $k$-subsets we need to consider.

We noted in Section 5.2 that for both parameter sets, $N_{\mathrm{Sym}(\mathcal{P})}(G) = H = \mathrm{AGL}(1, 41) \times \mathrm{AGL}(1, 11)$, with the action of $H$ on $\mathcal{P}$ as defined there. Further, $H$ is the largest possible subgroup of $\mathrm{Sym}(\mathcal{P})$ for which this is true. By Lemma 4.5.1 we can use the orbits of $H$ on 10-subsets as equivalence classes of 10-subsets of $\mathcal{P}$ without losing any information in our search. Recall $H$ consists of multiplications and additions on rows and columns.

We consider for now only parameter set 1. By Lemma 5.5.1, any starter block corresponds to mask 2A in Figure 5.4.1. This gives us a unique 2-row, both these points lying in 1-columns, and four 2-columns. Further, we know that rows of the 2 points in a given 2-column lie in one of 20 predetermined pairs. Suppose that $L$ is a starter block. We wish to find a new starter block, $\hat{L}$ that is in the same equivalence class as $L$ but has some extra known properties. Consider first the unique 2-row of $L$. We can map this to row 0 via some row addition in $H$. We now have a starter block, say $L_1$ which has row 0 as its unique 2-row. From here, we can map one of the points in the 2-row to column 0 via a column addition. This addition will preserve all rows, so the 2-row is still in row 0. We now have a starter block $L_2$, equivalent to $L$, containing the point $(0, 0)$ and some other point in row 0. Suppose the column of this other point is $f \neq 0$. Since $\mathbb{Z}_{11}$ is a field, $f^{-1} \in \mathbb{Z}_{11}$ will map $f$ to 1 (the inverse of $f$ in $\mathbb{Z}_{11}$). So there is a column multiplication in $H$ mapping this second point to the point $(0, 1)$. Furthermore, since multiplication fixes 0, this multiplication will preserve the point $(0, 0)$. This gives us a starter block $L_3$, equivalent to $L$ and containing $(0, 0)$ and $(0, 1)$. Now consider the 'leftmost' of the four 2-columns. That is, the 2-column with the lowest column number. This must be at least column 2 since we know that columns 0 and 1 are 1-columns. We cannot perform a column addition, or a column multiplication, since either of these would move one or both of $(0, 0)$ and $(0, 1)$. Similarly, we cannot perform a row addition, since both $(0, 0)$ and $(0, 1)$ would be moved. However, since row multiplications fix row 0, we can perform a row multiplication. We can thus map one of the points in this 2-column into row 1. Since the rows come in known pairs, this will also determine the row of the other point in this column to be $-1$. Call the result of applying these operations $\hat{L}$. Clearly $\hat{L}$ is equivalent to $L$, since we have performed 4 operations from $H$.

We have seen through this that given any starter block $L$, there is another starter block $\hat{L}$ that is equivalent to $L$, satisfying the following conditions:

- $(0,0) \in \hat{L}$
- $(0,1) \in \hat{L}$
- If $f$ is the lowest number 2-column in $\hat{L}$, then column $f$ contains the points in the rows corresponding to the row-pair $1, 40$.

It is thus sufficient to only consider potential starter blocks satisfying these conditions.

Given these conditions, we see there are very few choices to make for potential starter blocks in parameter set 1. First we choose 4 columns from $\{2, 3, \ldots, 10\}$ to be the 2-columns. The first of these has its rows decided. For the rest, we need to choose, in order, 3 row pairs of the remaining 19. This gives us $\binom{9}{4} \cdot 19 \cdot 18 \cdot 17 \approx 800000$ possibilities.

We now consider parameter set 2. We could consider each of the 7 different masks separately, performing a slightly different argument for each. Instead, we will give a general argument that can apply to any mask.

Again, take $L$ to be a potential starter block. We first consider $L$'s unique 2-row. As before, perform a row addition to make this 2-row row 0. This gives us starter block $L_1$, equivalent to $L$. In this case the columns containing the points from the 2-row might not be 1-columns. Of the two points, take the point whose column contains more points, or pick one arbitrarily if they contain the same number. Using a column addition, map this column to column 0. This gives us line $L_2$, which has row 0 as a 2-row and contains the point $(0,0)$. Now, as before, perform a column multiplication to map the second point in the 2-row to column 1. This gives us $L_3$, which contains the points $(0,0)$ and $(0,1)$ and column 0 contains at least as many points as column 1. Finally, take $f$ to be leftmost (lowest numbered) column containing a point from $L_3$, excluding columns 0 and 1. Take some point from column $f$ which is in $L_3$. Perform a row multiplication to map this point's row to row 1. This gives us $\hat{L}$, which satisfies the following conditions:

- $(0,0) \in \hat{L}$
- $(0,1) \in \hat{L}$
- Column 0 contains at least as many points as Column 1.
- Excluding the points $(0,0)$ and $(0,1)$, if $f$ is the numerically earliest column occupied by $\hat{L}$, then $(1,f) \in \hat{L}$.

We see that there are many more choices to be made for a potential starter block in parameter set 2.

## 5.7 Results of the Searches

We will not in this chapter consider the exact mechanism used in performing these searches. This work will be covered in Chapter 6. We will however consider the outcome of these searches here. For both parameter sets, no valid starter blocks were found. This means, assuming the search procedure was correct, an issue which we consider in Section 6.4, that there can be no 2-$(451, 10, 1)$ designs admitting either $(\mathbb{Z}_{41} \times \mathbb{Z}_{11}) \cdot \mathbb{Z}_5$ or $(D_{82} \times \mathbb{Z}_{11}) \cdot \mathbb{Z}_5$ as a line transitive and point imprimitive automorphism group, for any of the permutation groups of this type defined in Section 5.2.

## Conclusion

Here we have presented two possible parameter sets to search on. We have shown that group theoretic arguments may be used on very little information to construct all possible permutation groups given in the parameter sets. In both of these cases, we found there were four distinct permutation groups, so in both cases four separate searches need to be performed. In these cases however, the arguments to be made before searching were the same for all four groups. The only difference in the searches between the groups is their orbits on pairs, used for the Orbit Theorem requirements.

We have also introduced here the idea of a mask. Masks are useful when there are two or more different $G$-invariant partitions on the points. We first need to apply the Delandtsheer-Doyen Theorem (Theorem 4.4.2) and calculate the intercept vectors on each partition. After this has been done, the masks give an overall picture of the interaction between the number of points per class from each partition. Although, as we will see in 6.3, this information is not used in this search, it was particularly useful in our first parameter set, since we were able to show that only one mask of 7 could possibly result in a design.

Finally here we presented symmetry arguments for both the parameter sets. We saw that in both cases we were able to fix 2 of the 10 points in a starter block, and restrict the row location of 1 or 2 more points. This greatly reduces the number of sets we need to consider in the search.

As a result of this work, a search for designs conforming to these two parameter sets was run. In both cases the search returned a negative result, indicating that no $G$-point imprimitive and $G$-line transitive designs exist conforming to either parameter set.

# Chapter 6

# The Search Algorithm

## Introduction

Here we present the search algorithm that has been developed and used to search the second parameter set given in Chapter 5. The algorithm used to search the first parameter set is of less interest here, given the exceptional nature of the conditions on starter blocks. Using these conditions, a brute force search on 10-subsets satisfying the mask and row pair conditions outlined in Section 5.5 is sufficient.

We will outline here the two separate computational parts of the search. For the group theoretic computations GAP was used, since it is very powerful in these areas. The raw combinatorial computation was performed in C, since C is much faster than GAP for such work. Throughout this chapter, we will assume we are working with a given permutation group $G_i$ on $\mathcal{P}$ for some $i = 1, 2, 3, 4$. For each permutation group, the procedure is repeated. A full copy of the source code used for these computations may be found in Appendix A.

In Section 5.2 we labelled a point $(e, f)$ where $e$ is the point's row and $f$ its column. We numbered the rows and columns 0 to 40 and 0 to 10 respectively. When programming, it is more convenient to have points labelled by a single integer. This integer can then be used to reference arrays and lists. To this end, within GAP the point $(e, f)$ was numbered $f \times 41 + e + 1$ and within C it was numbered $f \times 41 + e$. The reason for this slight difference is that in GAP arrays are referenced to start at 1, in C they start at 0. When the C code reads points in, it subtracts 1 and when it outputs points it adds 1 to account for this. In this chapter we will consider the GAP numbering system, so points are numbered 1 to 451.

Notice this gives us an obvious ordering on points, which is equivalent to the column first then row ordering we considered in 5.2. We order the points based on increasing numbering. This leads to a natural, lexographic ordering of subsets of $\mathcal{P}$.

**Definition 6.0.1.** Let $A, B \subset \mathcal{P}$, where $A = \{a_1, a_2, \cdots, a_n\}$, and $B = \{b_1, b_2, \cdots, b_m\}$,

with $a_1 < a_2 < \cdots < a_n$ and $b_1 < b_2 < \cdots < b_m$. Then $A < B$ if there is some $i$ with $1 \leq i \leq min\{n, m\}$ such that $a_j = b_j$ for all $j = 1, 2, \cdots i - 1$ but $a_i < b_i$, or if $n < m$ and $a_j = b_j$ for all $j = 1, 2, \cdots n$.

## 6.1 Calculating Orbits On Pairs

As we saw in Section 4.3, the first step in using the Orbit Theorem to search a parameter set is to calculate the orbits of the given group $G_i$ on unordered pairs from the point set $\mathcal{P}$. Once a group has been constructed acting on a set, it is a simple matter to construct orbits of given points. This step is almost identical for both parameter sets. The only difference in the constructions is that $G_i$ in parameter set 1 is generated by an additional involution.

In Section 5.2 we found four generators for $G_i$ for parameter set 1. Recall these were $\alpha$, an addition by 1 to the row, $\beta$ an addition by 1 to the column, $a^{20}$, multiplication by $-1$ (modulo 41) on the row and $a^8 b^{2i}$, a multiplication of order 5 on the rows and the columns. These elements can all be constructed as permutations in GAP and hence the group $G_i$ generated by them constructed as a permutation group on $\mathbb{Z}_{41} \times \mathbb{Z}_{11}$ . The construction of $G_i$ for parameter set 2 is the same, except the generator $a^{20}$, the involution on rows, is omitted.

Once this group construction has been performed, the orbit containing each 2-subset of points is calculated. The 2-subsets are represented by $\{r, s\}$ where $1 \leq r < s \leq 451$. The set of all 2-subsets are traversed, one by one. When considering a 2-subset of points, if its orbit is already known, it is ignored. If the 2-subset's orbit is not known, it is calculated and all points in that orbit have their orbit recorded. As we saw in Section 5.5, there are 25 orbits on parameter set 1 and 45 on parameter set 2. Of these, 5 of the orbits on parameter set 1 contain 2255 points, the rest 4510 points. All the orbits on parameter set 2 contain 4510 points. The orbits are denoted by integers, 1 to 25 for parameter set 1 and 1 to 45 for parameter set 2.

This orbit information needs to be accessible by the combinatorial search. The orbit numbers are outputted to a file where the orbit of the 2-subset $\{r, s\}$ with $1 \leq r < s \leq 451$ is stored in line $r$, entry $s$.

## 6.2 Conditions on Starter Blocks

An important step in the application of the Orbit Theorem is working out how many 2-subsets from a starter block come from each orbit. At the present stage this calculation must be done by hand and the results inputted to the search algorithm. This calculation was discussed in Section 5.5.

There are two main conditions that must be satisfied by starter blocks. First is the orbit condition, given by the Orbit Theorem (Theorem 4.3.1). The second is the intercept vector condition given by the Delandtsheer-Doyen Theorem (Theorem 4.4.2). We base our search, as we will see in the next section, on checking $q$-subsets, with $q \leq k$. We thus require some similar conditions that can be applied to a $q$-subset of $\mathcal{P}$ in order for it to be contained in some starter block.

Suppose we have $n$ orbits on 2-subsets, say $O_1, O_2, \cdots, O_n$, and an orbit condition that a starter block must contain $b_i$ 2-subsets from $O_i$ for $i = 1, \ldots, n$, as given in Theorem 4.3.1. Suppose also the intercept vector condition on a given $G_i$-invariant partition $\mathcal{C}$ is $[d_0, d_1, \ldots, d_k]$. Recall by this we mean that for each starter block $L$ there are $d_i$ classes of $\mathcal{C}$ containing exactly $i$ points from $L$, for $i = 0, \ldots, k$.

From these two conditions we develop the partial orbit condition and the partial intercept condition.

**Definition 6.2.1.** A $q$-subset $A$, with $q \leq k$ is said to satisfy the *partial orbit condition* if $A$ contains no more than $b_i$ 2-subsets from the orbit $O_i$, for each $i = 1, \ldots, n$, and $A$ is said to satisfy the *partial intercept condition* if there are not more than $d_i$ classes of $\mathcal{C}$ containing exactly $i$ points from $L$, for $i = 0, \ldots, k$.

It is clear that if $A$ does not satisfy the two partial conditions then $A$ cannot be extended to a starter block by the addition of points. Also, any starter block must also satisfy the partial conditions. We refer to a $q$-subset satisfying the partial conditions as a *partial starter block.* Where we are dealing with $q$-subsets that may or may not satisfy the partial conditions, we will refer to them as *partial blocks.*

In this search we have two partitions, the rows and the columns. The intercept vector for the rows is $[32, 8, 1]$, giving us a unique 2-row. The symmetry conditions from Section 5.6 fix the two points in this 2-row. For adding all other points we need only ensure there is at most one point per row. There are two possible intercept vectors on the columns, $[4, 5, 1, 1]$ and $[5, 4, 2, 0]$. Checking the columns partial intercept condition is slightly more complex, so we introduce the column information of a $q$-subset. This information is useful in testing conditions on the columns.

**Definition 6.2.2.** Let $S$ be a $q$-subset of $\mathcal{P}$. The *column information* of $S$ is the tuple $\mathrm{ColInf}(S) = (c_0, c_1, \cdots, c_{10})$ where $\mathrm{ColInf}(S)_i = c_i = |S \cap (\text{column } i)|$.

This gives a way to check the partial column condition. We now require

$$|\{j \in \{1, 2, \cdots, 10\} | \mathrm{ColInf}(S)_j = i\}| \leq d_i \text{ for all } i \leq k.$$

## 6.3 The Algorithm

We consider now an algorithm to search for starter blocks on the point set $\mathcal{P}$. We note that the intercept vector condition is in fact two concurrent conditions for our search, since there is an intercept vector on the rows and one on the columns to consider.

The algorithm used here checks all possibilities for q-subsets of $\mathcal{P}$ with $q \leq k$ in a lexographic order. The ordering we use is as described in Definition 6.0.1.

The algorithm is defined recursively. At each stage in the algorithm, we have a partial starter block $A$. If $|A| = k$ then $A$ is a starter block, so the algorithm is complete and we can use $A$ to construct a point imprimitive, line transitive design. Otherwise, we will try to add a point to $A$ to create a new, larger partial starter block. There are several pieces of information we will need to do this.

Input to the algorithm will be the parameters of the design, $v$ and $k$, the column and row lengths, $c$ and $d$, the orbits on pairs information for the given group $G_i$, the intercept vectors on columns and rows, the initial partial starter block $I$ and any additional symmetry conditions. This initial partial starter block consists of all points that are known to be in any starter block, according to the symmetry conditions.

In this search, we have $I = \{1, 42\}$, corresponding to the points $(0, 0)$ and $(0, 1)$, as identified in Section 5.6. As noted earlier, these points make up the unique 2-row, so the row intercept vector is simple. The orbit information will be as calculated in GAP, as outlined in Section 6.1.

In order to avoid checking the same subset twice, we maintain, for each subset $A$, a set of points $Q$ which still need to be considered for addition to $A$. Once a point has been added to $A$ and tested, this point is removed from $Q$.

The recursive function is

```
Funct(A,Q)
    Input: A, Q ⊂ P with |A| < k.

    WHILE ( Q not empty )
        find a = min(Q)
        remove a from Q.
        IF ( A ∪ {a} violates the partial row, partial column or symmetry conditions )
            discard a
        ELSE
            Create B = A ∪ {a}
            # We must be careful not to alter the original A.
            IF ( B passes the partial orbit condition )
                IF ( |B| = k )
                    B is a starter block. Output B.
                ELSE
                    make P a copy of Q
                    # We need to copy Q so that subsequent calls to
                    # Funct do not alter it.
                    Funct(B,P)
        # If we reach this point, B could not be extended to a starter block.
        #Try adding another point to A from Q.
```

END WHILE
# Here $Q$ is empty, so $A$ cannot be extended to a starter block.
Output 'A cannot be extended to a starter block.'
END Funct

This algorithm will be proved correct in Lemma 6.3.2 First we will make a few comments about the notation used in the algorithm by running it on an example. Suppose $A$ is some set, $Q = \{1, 2\}$. We first consider appending 1 to $A$. If we were to just make the assignment $A = A \cup \{1\}$ then, when we consider appending 2 to $A$, we will get $A \cup \{1, 2\}$. This is the reason for using a new name for the partial block $B$.

Now suppose that $B = A \cup \{1\}$ passes all conditions. Suppose we were to run Funct($B, Q$). This call would remove 2 from $Q$ and try $B \cup \{2\}$. If this set was to fail, we would return to trying to add points to $A$. However, we have removed 2 from $Q$, so $Q$ is now empty. The algorithm would terminate without ever checking $A \cup \{2\}$. This is why we must make a copy of $Q$ for use in recursion.

### 6.3.1 Checking Conditions within the Algorithm

The heart of this function is testing if adding $a$ to $A$ would violate the row or column conditions. The row condition is easy to test; we need only make sure that $a$ does not belong to the row of any point in $A$. To check the column condition, we first find $j$, the column number of $a$ and $c = \text{ColInf}(A)$. We now construct a new column information, $c'$ by

$$c'_l = \begin{cases} c_l & \text{if } l \neq j \\ c_l + 1 & \text{if } l = j \end{cases} .$$

We can then check the column condition on $c'$, as described in Section 6.2.

Checking the symmetry condition will vary between searches. In this search, the only additional symmetry condition used is that the third point in $A$ (which will be the first point added to $I$) must come from row 1, as described in Section 5.6. We also saw in that argument that column 0 must contain more points than column 1. This condition was relaxed for the search as it would not have had a significant effect on the total time taken for the search, but it would increase its complexity.

We notice that the set of points to consider for addition to $B = A \cup \{a\}$ is copied from those of $A$. Doing this is valid. Clearly the set of points to consider for B must be a subset of the original $Q$. Furthermore, we need not consider any point from the original $Q$ which has already been removed. If there was some $a' \in Q$ that had been removed, then $A \cup \{a'\}$ cannot be extended to a starter block. Clearly then, $B \cup \{a'\} = A \cup \{a, a'\}$ cannot be extended to a starter block. Thus there is no need to consider $a'$ being appended to $B$.

The algorithm is started with Funct($I$,$\mathcal{P}\backslash$(row 0)). The initial $Q$ is chosen as $\mathcal{P}\backslash$(row 0) since we wish to consider consider all points from $\mathcal{P}$, but $I$ already contains two points from row 0 so no more points may be added from row 0.

**Lemma 6.3.2 (Correctness of Funct).** *Suppose that $A, Q \subset \mathcal{P}$ with $|A| < k$. Then:*

1. *The Funct$(A, Q)$, as described above, will terminate.*
2. *Funct$(A, Q)$ will output all starter blocks $L$ with $A \subset L \subset A \cup Q$.*
3. *If $L$ is outputted by Funct$(A, Q)$ then $L$ is a starter block and hence $(\mathcal{P}, L^{G_i})$ is a 2-$(v, k, 1)$ design that is $G_i$-point imprimitive and $G_i$-line transitive.*

**Proof.** 1. The algorithm terminating follows since at each recursive step we are reducing $|Q|$ by 1. Notice that $|Q|$ is initially finite and each call of Funct performs at most $|Q|$ operations.

2. Suppose the claim is false. Then there exist $A, Q \subset \mathcal{P}$ such that $A$ can be extended via $Q$ to a starter block $L$, but Funct$(A, Q)$ does not output $L$. Notice $|A| < k$ is bounded above. Let $k' < k$ be the maximum size of such an $A$. Take $A, Q \subset \mathcal{P}$, such that $|A| = k'$ and $L$ is an extension of $A$ via $Q$ with $L$ a starter block not outputted by Funct$(A, Q)$. Suppose $L = A \cup \{a_1, a_2, \cdots, a_m\}$, where $a_i \in Q$ for $i = 1, 2, \cdots, m$ and $a_i < a_j$ for all $i < j$. We will show that $L$ must be outputted, a contradiction. Notice that the while loop in Funct$(A, Q)$ considers each point from $Q$, since it only terminates when $Q$ is empty. Also, since $A \cup \{a_1\} \subseteq L$, it is a partial starter block. If $m = 1$ then $A \cup \{a_1\} = L$ and hence $L$ will be outputted. Suppose $m > 1$. Then the algorithm will call Funct$(B, P)$, where $B = A \cup \{a_1\}$ and $P \subset Q$. Since $a_2, a_3, \cdots, a_m > a_1$, we see that $\{a_2, a_3, \cdots, a_m\} \subset P$. This follows since the algorithm removes points from $Q$ in increasing order. Hence $L$ is an extension of $B$ via $P$. However, $|B| = k' + 1 > k'$, so $L$ must be outputted by the algorithm, by the maximality of $k'$. Hence, regardless of the value of $m$, the algorithm will always output $L$, a contradiction. Thus Funct$(A, Q)$ will output all starter blocks $L$ with $A \subset L \subset A \cup Q$, for all $A, Q \subset \mathcal{P}$.

3. If $L$ is outputted by the algorithm, then $|L| = k$ and $L$ is a partial starter block. Thus $L$ must be a starter block and hence $L$ satisfies the orbit condition. It follows, by the Orbit Theorem (Theorem 4.3.1), that $(\mathcal{P}, L^G)$ is a 2-$(v, k, 1)$ design. □

## 6.4 Testing the Search

Ensuring the search is running correctly is a difficult task, due to the large number of possibilities it must consider. There are three main areas where the search may fail due to human or computer failure. The search may not traverse the space of partial starter blocks correctly, and as a result fail to check some possibilities. Secondly, the conditions, either orbit or intercept, may pass false positive results, where a $q$-subset passes the condition incorrectly. Finally, a condition may give false negatives, where a partial starter block fails a condition when it should have passed.

In order to check these conditions, the partial blocks of size 6 were considered. Size 6 was chosen since, as we will see later, it is the largest size for which the number of partial blocks we need to consider is practical. If we were to try and consider partial blocks of size 7, any output would quickly become larger than any realistic storage media. It seems reasonable that if the algorithm is performing correctly on 6-subsets, it "should" be doing so at later stages in the search. Unfortunately, this is as rigorous as testing can be.

### 6.4.1 Testing the Intercept Vector Condition

The traversal of partial starter blocks and the partial intercept conditions are closely related. The algorithm is written so that no block is considered if it violates the partial intercept condition. In this case the partial intercept condition is broken into two conditions, a row and a column condition, which are checked independently.

Here we present the tests used when searching using the column intercept vector $[4, 5, 1, 1]$. The same tests we performed for the vector $[5, 2, 4, 0]$, although the number of partial blocks outputted at each stage varied.

**Lemma 6.4.1.1.** *The number of 6-subsets of $\mathcal{P}$ satisfying the partial intercept vector condition for vector $[4, 5, 1, 1]$ and the symmetry conditions is* 15719080*.*

**Proof.** The set of all possible column informations of a 6-subset conforming to the partial intercept vector conditions and the symmetry conditions was generated. In order for the symmetry conditions to be satisfied, we need at least one point each from columns 0 and 1. There were 380 such column informations. Given a column information, we need only select a row for each point to completely determine $S$. Of these 6 points, we know the rows of 3 of them by the symmetry conditions given in Section 5.6. We know the subset must contain the points $(0, 0)$ and $(0, 1)$ and that the next smallest point lies in row 1. This leaves 3 rows to determine. There are 3 possibilities for these points. They may all lie in the same row (i), 2 may lie in 1 row and 1 in a different row (ii) or all 3 may be in different rows (iii). Of the 380 column informations, 210 correspond to case (i), 145 to case (ii) and 25 to case (iii). For case (i) there are $\frac{39\times38\times37}{3!}$ possible choices for the rows, in case (ii) $\frac{39\times38\times37}{2!}$ choices and in case (iii) $39 \times 38 \times 37$ choices. This gives us $210 \times \frac{39\times38\times37}{3!} + 145 \times \frac{39\times38\times37}{2!} + 25 \times 39 \times 38 \times 37 = 15719080$ possibilities. □

When the algorithm is run without the partial orbit condition and told to output any 6-partial blocks it finds passing the partial intercept condition, it outputs exactly 15,719,080 partial blocks. This shows us that the algorithm is traversing correctly and that the partial intercept conditions are being correctly checked.

### 6.4.2 Testing the Orbit Condition

Testing the orbit condition relies on this output of all possible 6-partial blocks and the output of 6-blocks found that also satisfy the orbit condition. There were approximately 1,500,000 partial blocks found that passed the orbit condition according to the search, depending on the permutation group. The fact that these numbers were different for each group shows us these groups are indeed not isomorphic. First, these blocks were passed into GAP, where the orbit condition was rechecked using the orbit information calculated inside GAP, rather than the data outputted by GAP earlier. All of these partial blocks passed the orbit condition here, which tells us there were no false positives on the orbit condition.

Next, the set of 6-subsets that passed the partial intercept condition but not the orbit condition was generated. This was done by finding the difference between the output of all partial blocks and those found to pass the orbit condition. We ran through the list of partial starter blocks, removing each from the list of all possible 6-subsets. Since both lists are ordered, this can be done very quickly. All subsets remaining after this process had the partial orbit condition checked in GAP. None of these subsets passed this condition. This tells us there are no false negatives on the partial orbit condition.

This phase of testing verified that the algorithm was working correctly at a depth of 6. Although this does not tell us directly that the algorithm will function correctly up to a depth of 10, where it will terminate, it is good evidence that it will do so.

### 6.4.3 Testing the Algorithm on a Known Example

The final phase in testing was to test the algorithm on a different parameter set that has known results. We consider the Singer Cycle acting on the Desarguesian Projective Plane of order 7, as described in Section 3.5. Notice that $7^2 + 7 + 1 = 57 = 19 \times 3$, so we know there exists at least one such line transitive, point imprimitive, projective plane. In fact, it has been shown [1] that the Desarguesian Plane is the only projective plane of order 7 up to isomorphism. The corresponding parameter set is $v = 57, k = 7$ and $G = \mathbb{Z}_{19} \times \mathbb{Z}_3$ acting as addition on 3 rows and 19 columns. Although we do not consider the projective planes to be useful results in the overall search for point imprimitive, line transitive, linear spaces, this set is useful here. In particular, the small size of the point set makes the search very easy to run and to check. There were only two changes that needed to be made to the search code to run this search. First, the group $G$ needed to be constructed differently in GAP, but then the same GAP code was used for finding orbits on unordered pairs of points. Secondly, the constants for the number of rows and columns, as well as the length of a line and the number of orbits needed to be altered in the C code. This search was then run, it resulted in a single starter block being found, which generated a projective plane of order 7.

# Conclusion

The algorithm we have outlined here checks for starter blocks using a branch and cut traversal. Using this, it is possible to search a large number of possibilities (in this case of the order of 10 million) in a small time (here about 15 minutes). In most cases it is required to run the algorithm several times, over each permutation group ($G_i$) and each intercept vector on each $G_i$-invariant partition.

The algorithm used for this search is easily adaptable to a larger class of searches. First we require two orthogonal $G$-invariant partitions, treated as rows and columns as we have done here. If, in addition, one of the Delandtsheer-Doyen parameters $x$ and $y$ is 1 and the unique inner pair on the corresponding partition has been fixed by a symmetry argument, this algorithm can be used on the search. In order to use the algorithm here, we need only change some of the inputted constants, such as the number of orbits on 2-subsets, and the number of classes per partition.

# Chapter 7

# Conclusion

In this dissertation, we have examined the problem of identifying point imprimitive, line transitive, linear spaces. We have see that this is the same as identifying point imprimitive, line transitive, 2-$(v, k, 1)$ designs, which is in turn equivalent to searching for $k$-subsets of a set with cardinality $v$ satisfying certain properties.

Using the two main theorems considered here, the Orbit Theorem and the Delandtsheer-Doyen Theorem, in this manner is not new. These results were applied in the identification of almost all known examples of such designs which are not also projective planes. Using the Delandtsheer Doyen Theorem across multiple partitions to construct masks, which we introduced in Section 5.4 is a new approach to this problem. We saw here that in some cases, these masks can be used to greatly reduce the number of subsets that need to be considered.

Eight different groups acting on the same point set were examined. The search performed was for 2-$(v, k, 1)$ designs on this point set, admitting one of the eight possible permutation groups, as a line transitive, point imprimitive, automorphism group. It was shown that there are no such designs. As a result, these groups and the resultant parameter sets can be removed from consideration for future work.

Perhaps a more important result of this research is the algorithm used to perform the search on the last four groups. This algorithm can be used for a broader class of searches, any search where there are two orthogonal invariant partitions on the points and where $x$, the Delandtsheer Doyen parameter on one of those partitions, is 1.

The drive for future work in this area then is the identification of groups acting on some point set such that a 2-$(v, k, 1)$ design is likely to exist, admitting the group acting point imprimitively and line transitively.

# Appendix A

# Code for Searching

## A.1 Parameter Set 2 - $(\mathbb{Z}_{11} \times \mathbb{Z}_{41}) \cdot \mathbb{Z}_5$

The following is a listing of the C code used to perform the search for paramter set 2. The parameter set is described in Chapter 5, the algorithm used for the search is described in Chapter 6.

The code resides in three files: orbitwork.g, containing the GAP code to generate orbits on 2-subsets; design451.c, containing C code to perform the search; and datastructures.c, defining the datastructures used in the search. Detailed information about the internals of the data structures used here is not important, so only the header file datastructures.h is included.

## A.2 Running the Search

The two parts of the search were run on different machines. The GAP processing was performed in GAP4, running on an AMD K6-350, with 128MB of RAM. The C processing was run on an AMD K6-200, with 128MB of RAM. Both systems were running Red Hat Linux 7.0. For each intercept vector and permutation representation, the C based search took 15 minutes to run. Generating the orbit information in GAP took 10 minutes for each permutation representation.

## A.3 orbitwork.g

```
#This code sets up storage lists and functions for the calculation of the
#orbits of G=(Z41xZ11).Z5 on pairs of points taken from 1..451. It is based
#on code written by Alice Niemeyer.

npoints:=451;
f41 := GF(41);
f11 := GF(11);
y := PrimitiveRoot(f41);
z := PrimitiveRoot(f11);
el1 := Elements(f41);
el2 := Elements(f11);
omega:=[];
for i in [1 .. Length(el2)] do
    for j in [1 .. Length(el1)] do
        Add(omega, [el1[j], el2[i]]);
    od;
od;

perm1 := [];
for i in [1 .. Length(omega)] do
    perm1[i] := Position(omega, [omega[i][1]+One(f41), omega[i][2]]);
od;
p1 := PermList(perm1); #p1 is permutation from addition in f41

perm2:=[];
for i in [1 .. Length(omega)] do
    perm2[i] := Position(omega, [omega[i][1], omega[i][2]+One(f11)]);
od;
p2 := PermList(perm2); #p2 is permutation from addition in f11

a:=y^8;
b:=z^2;

permy:=[];
for i in [1 .. Length(omega)] do
    permy[i] := Position(omega, [omega[i][1]*a, omega[i][2]]);
od;
py:=PermList(permy);

permz:=[];
pz := [];
```

```
for i in [1 .. 4] do
    permz[i]:=[];
    for j in [1 .. Length(omega)] do
         permz[i][j] := Position(omega, [omega[j][1], omega[j][2]*(b^i)]);
    od;
    pz[i]:=PermList(permz[i]);
od;

G:=[];
for i in [1 .. 4] do
    G[i]:=Group(p1, p2, py*pz[i]);
od;

OrbitInfo:=[];
for i in [1 .. 4] do
    OrbitInfo[i]:=[];
    for j in [1 .. npoints] do
        OrbitInfo[i][j]:=[];
        for k in [1 .. (j-1)] do
            OrbitInfo[i][j][k]:=0;
        od;
    od;
od;

#Sets all the values in OrbitInfo and returns the number of orbits
#Returns the number of different orbits there are
setOrbitInfo := function(info, grp)
    local orbitNumber, orbit, i, j, p;
    orbitNumber:=1;
    for i in [1 .. npoints] do
        for j in [1 .. (i-1)] do
            if info[i][j]=0 then
                orbit := Orbit(grp, Set([i,j]), OnSets);
                Print("Length of orbit is ", Length(orbit), "\n");
                for p in orbit do
                    if p[1]<p[2] then
                        info[p[2]][p[1]] := orbitNumber;
                    else
                        info[p[1]][p[2]] := orbitNumber;
                    fi;
                od;
                Print("Done orbit number ", orbitNumber,"\n");
                orbitNumber:=orbitNumber+1;
```

```
            fi;
        od;
    od;
return orbitNumber-1;
end;

writeOrbitInfo := function(info, filename)
    local i,j;
    PrintTo(filename," ");
    for i in [2 .. Length(info)] do
        for j in [1 .. (i-1)] do
            AppendTo(filename,info[i][j]," ");
        od;
        AppendTo(filename,"\n");
    od;
end;

for i in [1..4] do
    setOrbitInfo(OrbitInfo[i],G[i]);
od;
writeOrbitInfo(OrbitInfo[1],"orbit.design.1");
writeOrbitInfo(OrbitInfo[2],"orbit.design.2");
writeOrbitInfo(OrbitInfo[3],"orbit.design.3");
writeOrbitInfo(OrbitInfo[4],"orbit.design.4");
```

## A.4 design451.h

```
#include <signal.h>
#include <stdio.h>
#include <string.h>
#include <time.h>
#include <malloc.h>

#define SN      451           /* the total number of points */
#define BS      10            /* the size of a block */
#define Ysize   41
#define Xsize   11
#define NO      45            /* the number of orbits on pairs */
#define TRUE     1
#define FALSE    0
#define MaxColLength 3

#include "datastructures.h"
```

## A.5 design451.c

```
/*
** Code to search for designs using 451 points with a 41 by 11 grid structure
** Based on origional code by Alice Niemeyer
** Greg Cresp
*/

#include "design451.h"
#include <time.h>
/* Some global variables */
unsigned long   OrbitNr[SN][SN];  /* orbits on 2-subsets */
unsigned long   Block[BS];        /* the block we hope to find */

/* Read the orbit on pairs, if (i,j) lies in k-th orbit set the
** k-th bit of OrbitNr[i][j]
** Code used unmodified from Alice Niemeyer
*/
void read_orbit(char *fileNumber) {
  unsigned int  i,j,n;
  FILE      *fp;
  char* filename;
  filename=(char*)malloc((strlen("orbit.design.")+1)*sizeof(char));
  strcpy(filename,"orbit.design.");
  strcat(filename, fileNumber);
  if( (fp = fopen( filename, "r" )) == NULL ) {
    perror( "orbit.design" );
    exit( 1 );
  }
  free(filename);
  for( i = 0; i < SN; i++ )
    for( j = 0; j < i; j++ ) {
      if( fscanf( fp, "%d", &n ) != 1 ) {
    fprintf( stderr, "reading orbit failed at %d %d\n", i, j );
    exit( 1 );
      }
      OrbitNr[i][j] = OrbitNr[j][i] = n-1; /* GAP starts at 1, C at 0 */
    }

  if( fscanf( fp, "%d", &n ) == 1 ) {
    fprintf( stderr, "number left : %d\n", n );
    exit( 1 );
  };
```

```
  fclose( fp );

  printf( "read orbit.design\n" );
  fflush( stdout );
}

/* Continue the search at level lev, using the covered information given.
   All columns up to and including start_col have already been considered.
   Pre: Data (global variable) up to lev still satisfies all conditions.
*/
void search( unsigned long lev, covered_orbits_data used,
    covered_rows_data used_rows,
    unsigned long current_col, column_length_data col_len,
    unsigned int length, current_length_data cur_len,
    unsigned int start_row) {

  int  row, i;
  covered_orbits_data used_copy;
  covered_rows_data used_rows_copy;
  column_length_data col_len_copy;

  /* if a block is full just print it */
  if( lev == BS ) {
    printf( "yeah { " );
    for( i = 0; i < BS-1; i++ )
      printf("%ld,", Block[i]+1 );
    printf("%ld ", Block[BS-1]+1 );
    printf( " },\n" );
    fflush( stdout );
    return;
  }

  /* This should only happen on the first column. It means we need to
     assign a length */
  if(length==-1) {
    /* get the lowest length. If this is -1 (it shouldn't ever be), stop */
    if((length=get_next_length(&col_len, get_current_length_data
                  (&cur_len,current_col)-1))==-1) {
      printf("Could not get new length\n"); fflush(stdout);
      return;
    }
    search(lev, used, used_rows, current_col, col_len, length,
       cur_len, 1);
```

```
}
/* if the current column is 'full' */
else if(length==get_current_length_data(&cur_len,current_col)) {
  /* See if we can extend the current column to a longer column. */
  col_len_copy=col_len;
  /* another column of the current length is allowed */
  col_len_copy.data[length]++;
  /* If we can find a larger length, continue */
  if((length=get_next_length(&col_len_copy, length))!=-1) {
    search(lev, used, used_rows, current_col, col_len_copy,
       length, cur_len, start_row);
  }
  /* Now try going to the next column */
  if(current_col<Xsize-1) {
    col_len_copy=col_len;
    /* get the lowest length. If this is -1 (indicates an invalid
   point arrangement), stop */
    if((length=get_next_length(&col_len_copy, get_current_length_data
               (&cur_len,current_col+1)-1))==-1) {
  printf("Could not get new length, %d\n",get_current_length_data
         (&cur_len,current_col+1) ); fflush(stdout);
    }
    else {
  search(lev,used,  used_rows, current_col+1, col_len_copy, length,
         cur_len, 1);
    }}}
/* if there are still points to add to the current column */
else if(length>get_current_length_data(&cur_len,current_col)){
  inc_current_length_data(&cur_len, current_col);
  /* The third point is always in row 1 */
  if(lev==2) {
    Block[lev]=current_col*Ysize+1;
    used_copy=used;
    used_rows_copy=used_rows;
    if(checkOrbits(lev, &used_copy)&&checkRows(1, &used_rows_copy)) {
  search(lev+1, used_copy, used_rows_copy, current_col, col_len,
         length, cur_len, 2);
    }
  }
  else {
    for(row=start_row;row<Ysize; row++) {
  Block[lev]=current_col*Ysize+row;
  used_copy=used;
```

```
    used_rows_copy=used_rows;
    if(checkOrbits(lev, &used_copy)&&checkRows(row, &used_rows_copy)) {
      search(lev+1, used_copy, used_rows_copy, current_col, col_len,
         length, cur_len, row+1);
}}}}}

/* add the first two points to the block and start the search */
void setUpBlockAndSearch(column_length_data col_len) {
  int lev = 0;
  int col = 0;
  covered_orbits_data used;
  covered_rows_data used_rows;
  current_length_data cur_len;
  /* make sure all the data structures are zeroed */
  clear_covered_orbits_data(&used);
  clear_covered_rows_data(&used_rows);
  clear_current_length_data(&cur_len);
  /* add the point (0,0) */
  Block[lev] = 0+col*Ysize;
  inc_current_length_data(&cur_len,col);
  lev++;
  col++;
  /* add the point in row 0, column 1 */
  Block[lev] = 0+col*Ysize;
  /*update the orbits information*/
  inc_current_length_data(&cur_len,col);
  checkOrbits(lev, &used);
  checkRows(0,&used_rows);
  lev++;
  /* do the search */
  search(lev,used,used_rows,0, col_len, -1, cur_len, 0);
  fflush(stdout);
}

/* initialise all the global variables and set up the block */
/* Program needs 4 arguments, <group num> <d_1> <d_2> <d_3> */
int main(int argc, char** args) {

  column_length_data col_len;
  int i,t,sum;
  if(argc!=5) {
    fprintf(stderr,"Insufficient arguments\n");
    exit(1);
```

```
  }
  if(strlen(args[1])!=1||args[1][0]<'1'||args[1][0]>'4') {
    fprintf(stderr,"Invalid arguments\n");
    exit(1);
  }
  clear_column_length_data(&col_len);
  sum=0;
  for (i=2;i<argc;i++) {
    col_len.data[i-1]=(int)args[i][0]-'0';
    sum+=(int)args[i][0]-'0';
  }
  /* fill in the required number of 0-columns */
  col_len.data[0]=Xsize-sum;
  init_globals();
  read_orbit(args[1]);
  t=(int)time(NULL);
  printf("Start time is %d\n",t);
  setUpBlockAndSearch(col_len);
  t=(int)time(NULL);
  printf("End time is %d\n",t);
  return 1;
}
```

## A.6 datastructures.h

```
#define LONG_INTS_FOR_USED_ORBITS (((NO+1) / (sizeof(long int) * 8)) + 1)
#define LONG_INTS_FOR_ROWS (((Ysize+1) / (sizeof(long int) * 8)) + 1)

struct covered_orbits_data {long int data[LONG_INTS_FOR_USED_ORBITS];};
struct covered_rows_data {long int data[LONG_INTS_FOR_ROWS];};
struct column_length_data {int data[MaxColLength+1];};
struct current_length_data {unsigned int data[Xsize];};

/* Stores the number of each length column still required */
typedef struct column_length_data column_length_data;
/* Stores the number of points currently in each column */
typedef struct current_length_data current_length_data;
/* Stores the orbits on pairs currently used */
typedef struct covered_orbits_data covered_orbits_data;
/* Stores the rows currently used */
typedef struct covered_rows_data covered_rows_data;

/* Initialises the global variables */
void init_globals();
/* Set or test the given bit of the orbit or row data  structure */
void orbs_set_bit(covered_orbits_data *covered, int bit);
int orbs_test_bit(covered_orbits_data *covered, int bit);
void rows_set_bit(covered_rows_data *covered, int bit);
int rows_test_bit(covered_rows_data *covered, int bit);
int checkOrbits(unsigned int level, covered_orbits_data *orbitsUsed);
int checkRows(unsigned long row, covered_rows_data *used_rows);
/* Zero the given data structure */
void clear_covered_orbits_data(covered_orbits_data *d);
void clear_covered_rows_data(covered_rows_data *d);
void clear_column_length_data(column_length_data *d);
void clear_current_length_data(current_length_data *d);
void inc_current_length_data(current_length_data *d, int col);
unsigned int get_current_length_data(current_length_data *d, int col);
/* Debug functions to print the given structure */
void printblock(int lev);
void printused(covered_orbits_data *used);
void printrows(covered_rows_data *used);
void print_current_length_data(current_length_data *col_len, int column);
void print_column_length_data(column_length_data *col_len);
/* Returns the first available column length greater than start */
unsigned int get_next_length(column_length_data *col_len, int start);
```

# Bibliography

[1] R. C. Bose and K. R. Nair, 'On complete sets of Latin squares', *Sankhyā* **5** (1941), 361–382.

[2] M. J. Colbourn and C. J. Colbourn, 'Cyclic Steiner systems having multiplier automorphisms', *Utilitas Math.* **17** (1980), 127–149.

[3] A. Delandtsheer and J. Doyen, 'Most block-transitive $t$-designs are point-primitive', *Geom. Dedicata* **29** (1989), no. 3, 307–310.

[4] J. D. Dixon and B. Mortimer, *Permutation groups* (Springer-Verlag, New York, 1996).

[5] D. R. Hughes and F. C. Piper, *Design theory* (Cambridge University Press, Cambridge, 1985).

[6] W. H. Mills, 'Two new block designs', *Utilitas Math.* **7** (1975), 73–75.

[7] W. Nickel, A. C. Niemeyer, C. M. O'Keefe, T. Penttila, and C. E. Praeger, 'The block-transitive, point-imprimitive 2-$(729, 8, 1)$ designs', *Appl. Algebra Engrg. Comm. Comput.* **3** (1992), no. 1, 47–61.

[8] W. Nickel and A. C. Niemeyer-Nickel, 'Two families of block-transitive, point-imprimitive block designs', *ANU Mathematical Sciences Research Centre Report* (1990), no. SMS-022-90.

[9] C. E. Praeger and A. C. Niemeyer, *Private Communication* (2001).

[10] C. E. Praeger and N. D. Tuan, 'Inequalities involving the Delandtsheer-Doyen parameters for finite line-transitive linear spaces', *J. Combin. Theory Ser. A* **102** (2003), 38–62.